\documentclass[12pt,a4paper]{article}
\usepackage[utf8]{inputenc}
\usepackage{amsmath}
\usepackage{listings}
\usepackage{graphicx}
\usepackage{float}
\usepackage{subfig}
\usepackage{geometry}
\usepackage{amssymb}  
\usepackage{cite} 
\usepackage[final]{hyperref}
\usepackage{xcolor}

\newtheorem{theorem}{Theorem}
\newtheorem{lemma}{Lemma}

\title{
Global stability of SIR model with heterogeneous transmission rate modeled by the Preisach operator}
\author{Ruofei Guan, Jana Kopfov\'a, Dmitrii Rachinskii}
\date{}

\begin{document}

\maketitle

\begin{abstract}
In recent years, classical epidemic models, which assume stationary behavior of individuals, have been extended to include 
an adaptive heterogeneous response of the population to the current state of the epidemic.  However, it is widely accepted that human behavior can exhibit history-dependence as a consequence of learned experiences.  This history-dependence is similar to hysteresis effects that have been well-studied in control theory. To illustrate the importance of history-dependence for epidemic theory,
we study dynamics of a variant of the SIRS model where individuals exhibit lazy-switch responses to prevalence dynamics.
The resulting model, which includes the Preisach hysteresis operator, possesses a continuum of endemic equilibrium
states characterized by different proportions of susceptible, infected and recovered populations. 
We discuss stability properties of the endemic equilibrium set and relate them to the degree of heterogeneity of the adaptive response. 
Our results support the argument that public health responses during the emergence of a new disease can have 
long-term consequences for subsequent management efforts.
The main mathematical contribution of this work is a method of global stability analysis, which uses a family of Lyapunov 
functions corresponding to different branches of the hysteresis operator.
\end{abstract}



\section{Introduction}

The standard SIR model assumes constant transmission rate (Kermack and McKendrick 1927). However, during epidemics, the transmission rate can change due to interventions of the health authorities 
and an adaptive response of the population to 
dynamics of the infectious disease. The role of the modulation of the transmission rate due to a human response to dynamics of epidemics has been widely studied (Capasso and Serio 1978; d'Onofrio and Manfredi 2009;
d'Onofrio and Manfredi 2020).
Preventive measures that decrease the transmission rate can include media campaigns
raising awareness in the population about the current severity of the epidemic, 
massive vaccination campaigns, 
access to 
effective and affordable 
 tests and medicines provided by the authorities, face covering in public places,
quarantine and social distancing measures, transition to online teaching at schools and colleges,  gathering and travel limitations, business restrictions
 and stay-at-home 
 orders (Gostin and Wiley 2020; Marquioni and de Aguiar 2020; Center of Disease Control and Prevention 2021). 
 Public participation in these measures can be partly mandatory and partly voluntary. 
 
One major objective of preventive measures 
is flattening the curve, 
i.e.\ slowing down the spread of 
infection in order to keep the number of active disease cases at a manageable level 
dictated by the  capacity of the health care system (Matrajt and Leung 2020). 
However, restrictive
preventive measures 
can be typically introduced only for a limited period of time 
due to the economic and social cost they incur (Fairlie 2020). 
Because of these constraints, many prevention policies (such as the intervention protocols which have been recently implemented by the health authorities in order to contain the Covid-19 epidemic (Wilder-Smith et al. 2020) can be thought of, at least simplistically, as 
threshold based.
As such, an intervention starts when a certain variable such as the number of daily new infections, the percentage of the occupied hospital beds or the effective reproduction number $R$ reaches a critical threshold value set by the health or government authority (Department of Health and Human Services, Nebraska 2020; DeBenedetto and Ruiz 2021); the intervention is revoked when this variable drops below the level deemed safer.
If the thresholds at which the intervention begins and ends are different, then {\em hysteresis} effect is present.

In their earliest form, hysteresis effects, which are commonly called a `backward bifurcation' (a subcritical steady state bifurcation) in the epidemic theory literature, were clearly identified by Hadeler and Chavez (1995) in their study of a core group. 
Hysteresis 
in adopting a `healthy behavior'
can also result from social contagion processes  
 such as social interaction, reinforcement and imitation (Su et al. 2017; Wang et al. 2017). 
For instance,
	in the presence of imperfect vaccine, 
hysteresis loops of vaccination rate arise with respect to dynamic changes in the perceived cost of vaccination (Chen and Fu 2019)
as
individuals revisit their vaccination decision in response to dynamics of the epidemic and through a social learning processes under peer influence. These hysteresis loops form a complex structure with minor loops nested within the maximal loop.

Hysteresis is a type of  dependence of the state of a system on its past states (the term was first coined to describe 
magnetic hysteresis as lagging of the magnetization of a ferromagnetic material behind variations of the magnetizing field (Ewing 1881)). Characteristically, 
this lagging
is determined by the past {\em values} of the phase variables and the {\em ordering} of the past events rather than the {\em timing} of those events. 
In particular, hysteresis is not associated with any explicit time scale and, as such, differs
from the delayed response of systems with deviating argument (delayed and functional differential equations). In fact, hysteresis is defined as a {\em rate-independent} operator in (Visintin 1994) and is characterized by 
a durable effect of temporary stimuli in biological systems and population processes (Pimenov et al. 2012).
\ Thus, the rate-independence of hysteresis and the prescribed time scale of systems with a deviating argument (and closely associated convolution operators) can be viewed as two idealizations which present a modeler with a choice of complementary mathematical tools for modeling a particular process. 

The effectiveness of the prevention measures depends on the response of the population to the intervention policies and, in particular, the degree of homogeneity of this response. 
The ability and willingness of an individual to receive immunization 
or to follow community isolation policies
depends on the level of trust to  the authorities in the community, personal beliefs, perceived risk of contracting the disease, risk of possible complications and other factors, which vary with age, health condition, living circumstances, profession and social group (Guidry et al. 2021; Lazarus et al. 2021).
All these variations 
contribute to the heterogeneity of the response of the society
to the advent of an epidemic
including behavioral changes.
Effects of switching self-protecting behavior on the transmission of a
disease were analyzed, for example, 
by Geoffard and Philipson (1996) and Chen (2004)  who modeled a population of agents with heterogeneous preferences for taking protective actions using rational choice theory.
Oscillatory vaccination dynamics in heterogeneous populations were identified by Reluga, Bauch and Galvani (2006).
An important finding 
by Chen and Fu (2019) who studied a relaxation oscillator limit of these dynamics is that {\em hysteresis in adopting healthy behavior} becomes more pronounced 
with increasing heterogeneity of the population.

In this paper, we attempt to model the effect of {\em hysteresis} and {\em heterogeneity} in the response 
of the population to preventive measures on the epidemic trajectory.
Our focus is on the global stability of the set of equilibrium states.
We use a variant of the SIR model where the transmission rate depends on dynamics of the infected population.
As a starting point, we adopt the approach used by Chladn\'a et al. (2020) to modeling the homogeneous switched response of the population to the varying number of infected individuals by a two-threshold two-state {\em relay operator}.
In this model, it is assumed that the health authorities implement a two-threshold intervention policy
whereby the intervention starts when the number of infected individuals exceeds a critical threshold value and stops whenever this number drops below a different (lower) threshold. 
The objective of the two-threshold policy is to navigate the system to the endemic equilibrium simultaneously keeping the number of infected individuals in check and not committing to continuous intervention. If the response is ideally homogeneous, then  prevention measures are assumed to
 translate immediately to reduced values of the transmission coefficient and the effective reproduction number during the intervention.
 This response is characterized by the simplest rectangular 
 hysteresis loop in the relationship between  the density of the infected population and 
 the transmission coefficient. The rectangular hysteresis loop can be considered as
 a model of the bi-stable response associated with the backward bifurcation.

Next, to reflect the heterogeneity of the response, 
the population is divided into multiple subpopulations,
each characterized by a different pair of switching thresholds. In order to keep the model relatively simple, we apply averaging under additional simplifying assumptions
as in Kopfov\'a et al. (2021); Rouf and Rachinskii (2021). 
This leads to a differential model with just two variables, $S$ and $I$, but with a complex operator relationship between the transmission rate and the density of the infected population $I$.  As such, this operator relationship, known as the Preisach hysteresis operator (see, for example, Krasnosel'skii and Pokrovskii 1989; Mayergoyz 2003) 
accounts for hysteresis in the heterogeneous response.
On the other hand, applications of differential models with a population of two-state relay operators to population dynamics date back to the pioneering work by J\"ager and Hoppensteadt (1980).

We show that due to hysteresis, the heterogeneous model has a {\em connected continuum} of endemic equilibrium states characterized by different proportions of the susceptible, infected and recovered populations. 
We consider the convergence of the epidemic trajectory to this continuum.
In mechanical systems, hysteresis incurs energy losses thus promoting convergence to an equilibrium state.
This convergence can be shown using the energy potential as a Lyapunov function (for the construction of hysteretic thermodynamic potentials and examples, see Brokate and Sprekels (1996); Krej\v{c}\'\i\ et al. (2013); Krej\v{c}\'\i\ et al. (2019);
related methods of the construction of the Lyapunov function for the problem of absolute stability of Lurie systems with hysteresis in the feedback loop can be found in
Gelig and Yakubovich (1978) and the review Leonov et al.(2017).
However, for SIR systems, hysteresis is a source of instability, which can lead to periodic oscillations corresponding to recurrent waves of infection (Chladn\'a et al. 2020; Kopfov\'a et al. 2021).
The  objective of this work is to show the convergence to the equilibrium set
if hysteresis has a limited magnitude (as quantified by an appropriate measure)
and provide estimates for the magnitude of the hysteresis effect, which guarantee this stable scenario globally.
A numerical case study of the dynamical scenarios exhibited by the heterogeneous SIR model considered below has been provided by Rachinskii and Rouf (2021).

To analyze stability, we 
interpret the SIR model with the Preisach operator as a switched system (di Bernardo et al. 2008) associated with an infinite
family of nonlinear planar vector fields. 
Between the switching moments, a trajectory of the system is an integral curve of a particular globally stable vector field with the associated global Lyapunov function. The Preisach operator imposes non-trivial rules for switching from one vector field  to another corresponding to switching from one to another branch
of the hysteresis nonlinearity.
We prove stability by carefully controlling (a)
the increments
of the Lyapunov functions along the trajectory between the switching points; and, (b) the difference 
of the Lyapunov functions corresponding to different vector fields (associated with different hysteresis branches) at 
the switching points. This methodology has been previously applied to a simpler special
situation where the switching surface 
was the same for all the hysteresis branches 
(Kopfov\'a et al. 2021\cite{Kopfova2021}). The main 
mathematical
contribution of this paper is the extension of this method of global stability analysis to the more typical 
generic case
when the switching surface changes dynamically depending on the history of the state.

The paper is organized as follows. In the next section, we discuss an SIR model with a hysteretic heterogeneous response.
Stability results are presented in Section 3 and discussed in Section 4. The Appendix contains the proofs.

\section{Model}

{\bf 2.1.} The systems considered below adapt the standard dimensionless SIR model
\begin{equation}\label{basemodel}
\begin{aligned}
    \dot I &= \beta SI - (\gamma + \mu)I,\\
    \dot S &= -\beta SI - \mu S + \mu
\end{aligned}
\end{equation}
with simple immigration demography,
where $I$ and $S$ are the densities of the infected and susceptible populations, respectively; $\beta$ is the transmission coefficient; $\gamma$ is the recovery rate; and $\mu$, is the departure rate
due to the disease unrelated death and emigration. We assume a constant total population scaled to unity, hence the density of the recovered population $R=1-I-S$ can be removed from the system. 
The domain
$I, S\ge 0$, $ I+S\le 1$ is positively invariant for system \eqref{basemodel}.

Re-scaling the time, one obtains the normalized system  
\begin{equation}\label{r0model}
\begin{aligned}
    \dot I &= R_0 SI - I, \\
    \dot S &= - R_0 SI - \rho S + \rho,
\end{aligned}
\end{equation}
with the dimensionless parameters 
$$
R_0 = \frac{\beta}{\gamma+\mu}, \qquad \rho = \frac{\mu}{\gamma+\mu}<1.
$$
Here $R_0$ is the
{\em basic reproduction number} and $R=R_0 S$ is the {\em effective reproduction number}
for any given density $S$ of the susceptible population. If $R_0<1$, then the {\em infection free} equilibrium $(I^*,S^*)=(0,1)$
is globally stable in the closed positive quadrant. On the other hand, if $R_0>1$, then
the infection free equilibrium is a saddle, and the positive {\em endemic} equilibrium $(I_*,S_*)$ defined by
\begin{equation}\label{equil}
I_* = \left(1-\frac{1}{R_0}\right)\rho, 
\qquad S_* = \frac{1}{R_0}
\end{equation}
is globally stable in the open positive quadrant.
If 
\begin{equation}\label{focus}
\rho(R_0)^2 < 4(R_0 - 1), 
\end{equation}
then the endemic equilibrium is of focus type.

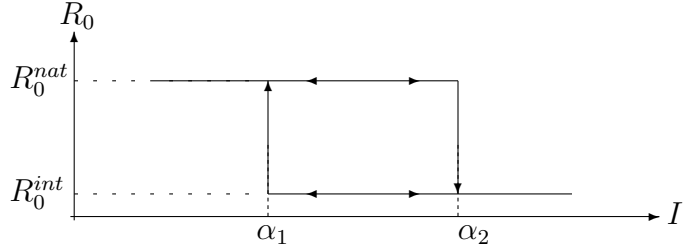
\begin{figure}[h]
\begin{center}
\setlength{\unitlength}{0.5mm}
\begin{picture}(120,30)
\put(0,30){\line(1,0){40}}
\put(55,30){\vector(-1,0){15}}
\put(55,30){\vector(1,0){15}}
\put(70,30){\line(1,0){10}}
\put(80,30){\vector(0,-1){30}}
\put(110,0){\line(-1,0){40}}
\put(30,0){\line(1,0){10}}
\put(55,0){\vector(-1,0){15}}
\put(55,0){\vector(1,0){15}}
\put(-21,-7){\vector(0,0){50}}
\put(-25,45){$R_0$}
\put(-22,-6){\vector(1,0){155}}
\put(135,-8){$I$}
\put(30,0){\vector(0,1){30}}
\put(-38,28){$R_0^{nat}$}
\multiput(-21,30)(5,0){10}{\line(1,0){1}} 
\put(-38,-2){$R_0^{int}$}
\multiput(-21,0)(5,0){10}{\line(1,0){1}} 
\put(27,-12){$\alpha_1$}
\multiput(30,-6)(0,2){10}{\line(0,1){1}} 
\put(80,-12){$\alpha_2$}
\multiput(80,-6)(0,2){10}{\line(0,1){1}} 
\thicklines
\end{picture}
\end{center}
\caption{{Switching rule \eqref{rev} for the basic reproduction number. The vertical segments correspond to instantaneous transitions of the basic reproduction number from the value $R_0^{nat}$ to $R_0^{int}$ and back.}
\label{relay1fig}}
\end{figure}

\medskip
\noindent
{\bf 2.2.} A switched system with two flows 
of the form (\ref{basemodel}), which have different values of the basic reproduction number,
$R_0=R_0^{nat}$ and $R_0=R_0^{int}$, was considered in Chladn\'a et al., 2020; Rouf and Rachinskii, 2021. Following this work, we assume that 
the parameter $\rho$ is the same for both flows 
and a switch from one flow to the other occurs when the density of the infected population reaches certain thresholds, $I=\alpha_1$ and $I=\alpha_2$,
satisfying 
$0<\alpha_1 < \alpha_2<1. $
More precisely, it is postulated that the basic reproduction number instantaneously switches from 
the value $R_0^{nat}$ to the value $R_0^{int}$ when the variable $I=I(t)$ reaches the upper threshold value $\alpha_2$.
On the other hand, $R_0$ switches back from the value $R_0^{int}$ to the value $R_0^{nat}$ as $I(t)$ reaches the lower
threshold value $\alpha_1$,  see Figure \ref{relay1fig}.
%
%
%

\bigskip
\bigskip
\bigskip



This two-threshold switching rule can be 
formalized 
%
%
using the standard {\em non-ideal relay operator}, which is
 also known as {a rectangular hysteresis loop} or a lazy switch (see e.g.~Visintin 1994).
 The relay is characterized by the threshold parameters $\alpha_1,\alpha_2$,
 and we use the notation $\alpha=(\alpha_1,\alpha_2)$.
The function $I(t): \mathbb{R}_+\to\mathbb{R}$ is called the {\em input} of the relay. The {\em state} of the relay, denoted $r_{\alpha}(t)$, equals either $0$ or $1$ at any moment $t\in \mathbb{R}_+$.
 Specifically,
 given any continuous input $I(t): \mathbb{R}_+\to \mathbb{R}$ and an initial value of the state, $r_{\alpha}(0)=r_{\alpha}^0$,
which satisfies the constraints
 \begin{equation}\label{v1}
 r_{\alpha}^0\in\{0,1\} \quad \text{if} \quad \alpha_1<I(0)< \alpha_2;
 \end{equation}
 \begin{equation}\label{v2} 
 r_{\alpha}^0=0 \quad \text{if} \quad I(0)\le \alpha_1;\qquad r_{\alpha}^0=1 \quad \text{if} \quad I(0)\ge\alpha_2,
 \end{equation}
 the state of the relay at the future moments $t>0$ is defined by the relations
 \begin{equation}\label{relay'}
 r_{\alpha}(t)=
 \left\{
 \begin{array}{cl}
 0 & \text{if there is $t_1\in [0,t]$ such that $I(t_1)\le \alpha_1$}\\
 & \text{and $I(\tau)<\alpha_2$ for all $\tau\in (t_1, t]$;}\\
 1 & \text{if there is $t_1\in [0,t]$ such that $I(t_1)\ge \alpha_2$}\\
 & \text{and $I(\tau)>\alpha_1$ for all $\tau\in (t_1, t]$};\\
 r_{\alpha}(0) & \text{if $\alpha_1<I(\tau)<\alpha_2$ for all $\tau\in[0,t]$.}
 \end{array}\right.
 \end{equation}
Function \eqref{relay'} 
will be denoted by
 \begin{equation}\label{re}
 r_{\alpha}(t)=({\mathcal R}_{\alpha}[r_{\alpha}^0]I)(t),\qquad t\geq 0.
 \end{equation}
 This standard short notation adopted from Brokate and Sprekels, 1996
 stresses that the relay's state is dependent on the history of $I(t)$, i.e.\ \eqref{relay'} defines an operator which takes the function $I(t): \mathbb{R}_+\to \mathbb{R}$ rather than an instantaneous value of $I$ as an input.
 Moreover,  the function $r_\alpha(t): \mathbb{R}_+\to \{0,1\}$ depends both on the input $I(t): \mathbb{R}_+\to\mathbb{R}$ and the initial state  $r_{\alpha}(0)=r_{\alpha}^0$  of the relay.
 By definition \eqref{relay'} of the input-to-state map \eqref{re}, the state satisfies the constraints
 \begin{equation}\label{compatibility}
 r_{\alpha}(t)=1 \quad \text{whenever} \quad I(t)\geq
 \alpha_2; \quad r_{\alpha}(t)=0 \quad \text{whenever} \quad I(t)\leq \alpha_1
 \end{equation}
 at all times. Further,
 the function \eqref{relay'} has at most a finite number of jumps between the values $0$ and $1$
 on any finite time interval $t_{0}\leq t\leq t_{1}$. 
 
 

Using formula \eqref{relay'} and notation \eqref{re}, 
the switching rule for the basic reproduction number depicted in Figure \ref{relay1fig} can be expressed 
as
\begin{equation}\label{rev}
R_0(t) = R_0^{nat} - (R_0^{nat}-R_0^{int})\cdot (\mathcal{R}_{\alpha}[r_{\alpha}^0]I)(t),\qquad t\ge0,
\end{equation}
Combining equations \eqref{r0model} with formula \eqref{rev} for $R_0=R_0(t)$,
one obtains a switched system. 
According to the interpretation discussed in the Introduction, the switched system \eqref{r0model}, \eqref{rev} 
models the dynamics of the epidemic under the assumption that the health authorities implement the two-threshold intervention policy --- 
an intervention begins when the density of infection (the number of active cases) exceeds the threshold value $\alpha_2$ 
and is revoked when the density of infection drops 
below the lower 
threshold value $\alpha_1$.
It is assumed that the intervention quickly translates into the reduction of the basic reproduction number; 
once the intervention stops,
$R_0$ returns to the larger value. 

As shown in Chladn\'a et al. (2020); Rouf and Rachinskii (2021),
each trajectory of switched system 
\eqref{r0model}, \eqref{rev} from the invariant domain 
$I, S\ge 0$, $ I+S\le 1$ converges to either an equilibrium state or a periodic orbit.
In particular, depending on the positioning of the thresholds $\alpha_1, \alpha_2$,
the system can have either one attractor (an equilibrium state or a periodic orbit), or
two attractors (two co-existing stable equilibrium states or a stable equilibrium state co-existing with a stable periodic orbit),
or three attractors (two equilibrium states and one periodic orbit).
For a periodic orbit with $I(t)$ oscillating between the values $\alpha_1$ and $\alpha_2$,
the point $(I(t),R_0(t))$ moves clockwise along the {rectangular hysteresis loop} shown in Figure~\ref{relay1fig}. 

\medskip\noindent
{\bf 2.3.} Model \eqref{rev} corresponds to an ideally homogeneous response of the population
when all the individuals change their behavior simultaneously following the recommendation of the health authority.
Now, we consider a model, in which several responses of the form \eqref{rev}, with different thresholds 
$\alpha$, are combined 
because different groups of individuals respond differently to the advent and dynamics of an epidemic
%
and to the interventions of the health authorities.

In particular, the ability and willingness to 
follow the recommendations of the health authority
can vary significantly from one to another group of individuals for the same level of threat of contracting the disease.
Multiple factors are at play such as the occupation, age, living environment and health conditions of an individual, to mention a few.
In order to account for the heterogeneity of the individual response, let us divide the susceptible population into non-intersecting subpopulations $S_\alpha$ parameterized by points $\alpha$ of a subset $\Pi\subset \{\alpha=(\alpha_1,\alpha_2): \alpha_1<\alpha_2\}$ of the $\alpha$-plane, assuming a homogeneous response within each group $S_\alpha$. As a simplification, let us assume 
that the response of the infected population 
is homogeneous, 
i.e.\ all infected individuals act the same and don't change their behaviors
over time,
hence the heterogeneity of the transmission coefficient is due to the behavior of the susceptible individuals only. 
Further, assume that 
the basic reproduction number for the subpopulation $S_\alpha$ is given by \eqref{rev}
and the contact structure is separable. 
Then, the average basic reproduction number for the entire population at time $t$ equals
\begin{equation}\label{rev''}
R_0(t)= R_0^{nat} 
-(R_0^{nat}-R_0^{int})\iint_{\Pi} 
({\mathcal R}_{\alpha}[r_{\alpha}^0]I)(t)\, d F(\alpha), \qquad t\in\mathbb{R}_+,
\end{equation}
where the 
cumulative distribution function $F(\alpha)$ describes the distribution of the susceptible {population} $S$ over the index set $\Pi$ (the set of threshold pairs $\alpha=(\alpha_1,\alpha_2)$). 
The effective reproduction number at time $t$ is then $R(t)=R_0(t)S(t)$.
Finally, we assume for simplicity that the distribution function $F(\alpha)$ is independent of time,
i.e.~$F(\alpha)$ does not change with variations of $I$. In this case, 
the mapping of the space of continuous
inputs $I(t):\mathbb{R}_+\to \mathbb{R}$ to the space of outputs $R_0(t):\mathbb{R}_+\to \mathbb{R}$
defined by \eqref{rev''} is known as the Preisach operator (see e.g.~Krasnosel'skii and Pokrovskii 1989).

Under the above assumptions, the dynamics of the epidemic is modeled by system \eqref{r0model}
where $R_0(t)$ is related to $I(t)$ by the Preisach operator \eqref{rev''}, which accounts 
for the heterogeneity of the transmission coefficient among different subgroups of the susceptible population. 

A similar system results from the assumption 
that the health authorities 
have multiple intervention policies (numbered $n=1,\ldots,N$) 
in place, each decreasing the transmission coefficient by a certain amount $ \Delta \beta_n$ while the intervention is implemented
(Rouf and Rachinskii 2021).
The authorities aim to provide an adaptive response, which is adequate to the severity of the epidemic.
Let us suppose that each  intervention policy 
is guided by the two-threshold start/stop rule, such as in \eqref{rev}, associated with a particular pair of thresholds $\alpha^n=(\alpha_1^n,\alpha_2^n)$. 
The response of the population is assumed to be homogeneous.
Under these assumptions, the basic reproduction number of system \eqref{r0model} is given by
\begin{equation}\label{rev'}
R_0(t)= R_0^{nat} 
- \sum_{n=0}^{N-1} \hat q_n\cdot
({\mathcal R}_{\alpha^n}[r_{\alpha^n}^0]I)(t),
\end{equation}
where $\hat q_n=\Delta\beta_n/(\gamma+\mu)$;
hence, $R_0$ is set to change at multiple thresholds $\alpha_1^n, \alpha_2^n$.
Operator \eqref{rev'}, known as the discrete Preisach model, is a particular case of \eqref{rev''} with an atomic measure $F$. 

Below we consider absolutely continuous distributions $F(\alpha)$ (see Appendix 5.1 for details). 
The corresponding operator \eqref{rev''}, which is called the continuous Preisach model,
can be written in the equivalent form
\begin{equation}\label{pre}
R_0(t)=R_0^{nat} - (R_0^{nat}-R_0^{int})\iint_{\Pi} q(\alpha)\, \bigl({\mathcal R}_{\alpha}[r_{\alpha}^0]I\bigr)(t) \,d\alpha_1 d\alpha_2,\qquad t\ge 0,
\end{equation}
where $0<R_0^{int}<R_0^{nat}$ and $q=q(\alpha): \Pi\to \mathbb{R}_+$ is a strictly positive measurable and bounded probability density function, i.e.
\begin{equation}\label{q}
\iint_\Pi q(\alpha)\,d\alpha_1d\alpha_2=1.
\end{equation}
This operator can be approximated by discrete operators \eqref{rev'}. We will assume that $q$ is bounded.
Also, in what follows,
\begin{equation}\label{pi}
\Pi=\{\alpha=(\alpha_1,\alpha_2): 0\le \alpha_1<\alpha_2\le1\}
\end{equation}
and
\begin{equation}\label{RR'}
R_0^{nat}>R_0^{int}>1.
\end{equation}

\begin{figure}[H]
    \centering 
    \includegraphics[width=0.8\textwidth]{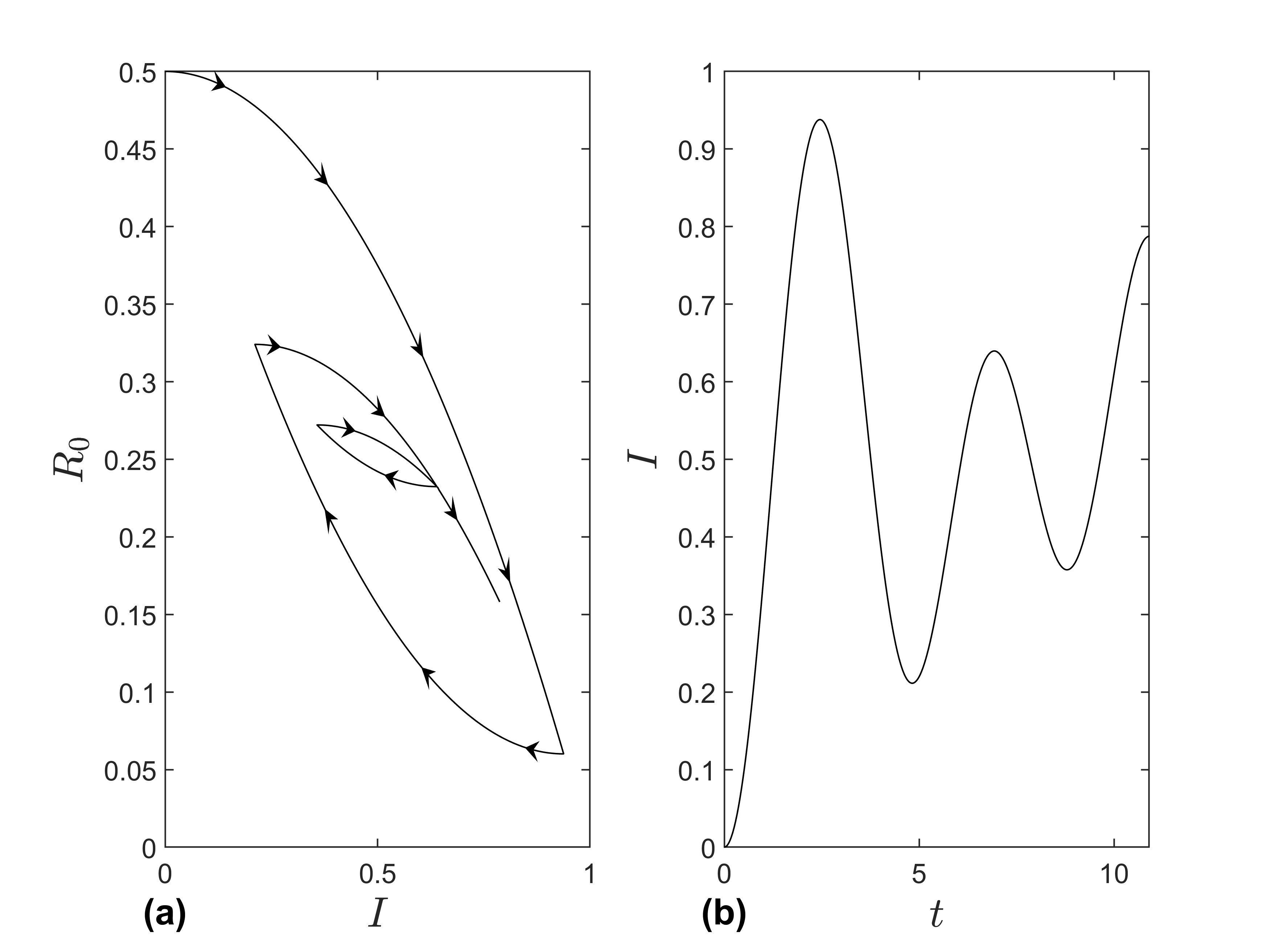} 
    \caption{An example of the input-output diagram for the Preisach operator \eqref{pre} on panel (a) for the input $I(t)$ shown on panel (b). The turning points on panel (a) correspond to local extrema of the input $I(t)$ on panel (b). \label{fig0}}
    \end{figure}

The input-output diagram of the Preisach operator in Figure \ref{fig0} illustrates the relation between $I(t)$ and $R_0(t)$
defined by equation \eqref{pre}. One can see that the point $(I,R_0)$ switches from one {\em branch} $R_0(I)$
to another at every turning point of the input $I(t)$. The presence of multiple branches and hysteresis loops
is a manifestation of the history-dependence, i.e.\ the output value $R_0(t)$  at a time $t$ is determined not only by the concurrent value $I(t)$ of the input but also by some of the past input values $I(\tau)$ achieved on the time interval $\tau\in [0,t]$. Fundamental properties of the history-dependence, which are typical of hysteresis operators, are briefly discussed in Appendix 5.1.


\section{Results}

{\bf 3.1.} Trajectories of system \eqref{r0model} coupled with the operator relationship \eqref{pre} lie in the infinite-dimensional phase space $\mathfrak U$ of triplets $(I,S,r_\alpha)$, where
the measurable function $r_\alpha=r(\alpha):\Pi\to \{0,1\}$ of the variable $\alpha=(\alpha_1,\alpha_2)$
describes
the states of the relays ${\mathcal R}_{\alpha}$ at a given moment 
(see Appendix 5.2 for details).
Slightly abusing the notation, we will 
also refer to the two-dimensional curve  $(I(t),S(t))$ as a trajectory, omitting the component \eqref{re} in the state space of the Preisach operator.

Let us first consider equilibrium states of system \eqref{r0model}, \eqref{pre}. The components $I$, $S$, $r_\alpha$ of the solution and the basic reproduction number $R_0$ at an equilibrium state are constant, and $R_0$ is related to the function $r_\alpha=r(\alpha):\Pi\to \{0,1\}$ by the equation
\begin{equation}\label{preR0}
R_0=R_0^{nat} - (R_0^{nat}-R_0^{int})\iint_{\Pi} q(\alpha) r_{\alpha} \,d\alpha_1 d\alpha_2.
\end{equation}
Due to assumption \eqref{pi} and the compatibility constraint \eqref{compatibility}, the inclusion $\alpha\in\Pi$ implies that all the relays are in state $r_\alpha=0$ when $I=0$. Therefore,
system \eqref{r0model}, \eqref{pre} has a 
unique {\em infection free} equilibrium state 
\begin{equation}\label{infectionfree}
E^*=(I^*,S^*)=(0,1),
\end{equation}
in which the function $r_\alpha: \Pi\to\{0,1\}$ is the identical zero and the basic reproduction number equals $R_0=R_0^{nat}$ according to \eqref{preR0}.
On the other hand, it is easy to see that due to assumption \eqref{RR'},
system \eqref{r0model}, \eqref{pre} also has a connected {\em continuum} of {\em endemic} equilibrium states
\begin{equation}\label{endemic}
E^{\theta}_*=(I_*^{\theta},S_{*}^{\theta})=\left(\left(1-\frac{1}{R_0^{\theta}}\right)\rho,\frac1{R_0^{\theta}}\right), \qquad 0\le\theta\le 1,
\end{equation}
which are characterized by different proportions of the infected, susceptible and recovered populations
and different values $R_0^\theta$ of the basic reproduction number (Rachinskii and Rouf 2021). 

\medskip
\noindent
{\bf 3.2.} The phase space $\mathfrak U$ 
	of system \eqref{r0model}, \eqref{pre} is naturally embedded into the space $\mathbb{R}^2\times L_1(\Pi;\mathbb{R})$
	and is endowed with a metric by this embedding (see Appendix 5.2 for details). As such, local and global stability properties of trajectories of system \eqref{r0model}, \eqref{pre}
	are defined in the standard way. We consider the positively invariant set 
	\[
	\mathfrak{U}_0=\{(I,S,r_\alpha)\in \mathfrak{U}: I>0, S>0, I+S\le 1\}.
	\]
	
	\begin{theorem}\label{t1}
For any given $\rho>0$, $R_0^{nat}>1$ and a strictly positive probability density function $q=q(\alpha): \Pi\to \mathbb{R}_+$ of the Preisach operator,
there is an $\varepsilon\in (0, R_0^{nat}-1)$ such that if $1+\varepsilon \le  R_0^{int} < R_0^{nat}$, then
	the set of endemic equilibrium states of system \eqref{r0model}, \eqref{pre} is globally asymptotically stable in $\mathfrak{U}_0$.
	Moreover, each trajectory from $\mathfrak{U}_0$ converges to an endemic equilibrium state.
	In particular, the components $I(t),S(t)$ of any trajectory from $\mathfrak{U}_0$ satisfy
	$(I(t), S(t))\to E_*^\theta$ for some $\theta\in[0,1]$.
	\end{theorem}
	
    


The proof presented in the Appendix provides an explicit estimate for the interval of values of $R_0^{int}\in [1+\varepsilon , R_0^{nat})$, for which the set of endemic equilibrium states is globally stable. In particular, the global stability is guaranteed if
 the quantities $\epsilon_0$ and $\kappa$ defined by \eqref{eps0}, \eqref{kappa} using \eqref{K}, \eqref{estim} are both positive.
	
\section{Conclusion}
We considered an SIR model where the transmission coefficient changes in response to
dynamics of the epidemic. 
We aimed to emphasize that the adaptive response of the population to the state of the epidemic is typically heterogeneous and that this response can feature history-dependence in a form, which is similar to a hysteresis effect. 
Specifically, we adapted a two-state two-threshold hysteretic switch as a model of
the adaptive response of an individual to the varying number of active cases.
In an ideally homogeneous population, the two switching thresholds of the transmission coefficient
can be imposed by 
the health authority which starts the intervention when the number of active cases exceeds a threshold $\alpha_1$ and ends the intervention when the number of active cases drops below another threshold $\alpha_2$. 
In order to account for the possibility of a heterogeneous 
response among the susceptible individuals,
we allowed a distribution of switching thresholds and
modeled 
the aggregate response of the susceptible population by the Preisach operator.
The variance of the distribution, $\sigma^2$,
measured the degree of heterogeneity of the public response 
to the interventions. 

The resulting heterogeneous model 
is shown to have 
a continuum of endemic equilibrium states 
differing by the proportions of susceptible, infected and recovered populations.
Numerical simulations of the model provided earlier by Rouf and Rachinskii, 2021 suggested that 
increasing the heterogeneity of the response
by widening the
spread of thresholds of different population groups (a larger variance $\sigma^2$) leads to a significant increase of the peak of infection during the epidemic. Moreover, a higher degree of heterogeneity of the public response also tends to steer the epidemic trajectory to an endemic equilibrium state with a higher 
prevalence and a lower proportion of the susceptible population. In other words, a more homogeneous response of the public to transmission prevention measures helps ``flattening the curve" and can lead to lower prevalence and transmission rate 
when an endemic equilibrium state is reached after the epidemic. These results suggest that 
intervention programs are more effective when accompanied by education campaigns, 
which convey the importance of the intervention measures to the public thus ensuring a more homogeneous response from different population groups.

However, when the history-dependent two-threshold response in our model is ideally homogeneous, i.e.\ the transmission coefficient behaves as a lazy-switch,
some threshold pairs lead to the convergence of the epidemic trajectory to a periodic orbit
predicting  recurrent outbreaks of the epidemic 
(Chladn\'a et al. 2020).
Emergence of sustained periodic oscillations of the prevalence 
was also shown numerically in a similar SIR model with vaccination 
when the heterogeneity of the response 
(measured by the variance $\sigma^2$ of the probability density function of the Preisach operator) decreases (Kopfová et al. 2021).
This scenario is in line with earlier findings by Reluga, Bauch and Galvani (2006)
who identified stable oscillatory vaccination dynamics in heterogeneous populations
and concluded that
the more homogeneous the response of a population, the less likely vaccine uptake is to be stable. 
The stability theorem presented in this paper agrees with all these results suggesting 
that the heterogeneity of the adaptive response promotes the convergence of the epidemic trajectory to an equilibrium state.
Indeed, the conditions of the theorem, which ensure the global
stability of the endemic equilibrium set, require a sufficient degree of heterogeneity of the response.
On the other hand, for a sufficiently homogeneous response, these conditions are violated leaving a possibility
of undesirable non-equilibrium dynamics,
i.e.\ a too homogeneous response can potentially trigger epidemic instability.

Thus, our results support the argument that a more homogeneous public response, which helps ``flattening the curve''
and tends to lead to a lower prevalence at the endemic state after the epidemic, comes at the price of a slower transitioning
to the endemic state and a higher risk of epidemic instability.
A higher degree of heterogeneity of the public response can stabilize and accelerate 
transitioning to an endemic state but 
after a higher infection peak; it also potentially results in a higher prevalence at the endemic equilibrium state
after the epidemic.

The mathematical contribution of this paper is a method of global stability analysis,
which uses a family of Lyapunov functions corresponding to different branches of the hysteresis operator.
We applied this method for the first time to an epidemic model with an adaptive transmission rate.
The method is potentially useful for analysis of other systems with hysteresis operators,
which have a continuum of equilibrium states.
  

\section{Appendix}

{\bf 5.1.~Continuous Preisach model.}
Let us briefly recall 
a rigorous definition of the 
continuous Preisach operator \eqref{pre} (Krasnosel'skii and Pokrovskii 1989). 
It involves a collection of non-ideal relays ${\mathcal R}_{\alpha}$, which respond to the same  continuous input $I(t): \mathbb{R}_+\to \mathbb{R}$ independently according to formula \eqref{relay'}. The relays contributing to the system have different pairs of thresholds $\alpha=(\alpha_1,\alpha_2)\in \Pi$,
where the subset $\Pi$ of the half-plane $\{\alpha=(\alpha_1,\alpha_2): \alpha_1<\alpha_2\}$ 
is assumed to be measurable and bounded;
the $\alpha$-plane is called the Preisach plane.
The output 
of the continuous Preisach model is the scalar-valued function $R_0(t): \mathbb{R}_+\to \mathbb{R}$
defined by \eqref{pre}, 
where 
$q(\alpha):\Pi\to\mathbb{R}_+$  is a positive bounded measurable function (probability density function) representing the weights of the relays; and, $r_{\alpha}^0$
is the initial state of the relay ${\mathcal R}_{\alpha}$ for any given $\alpha\in\Pi$.  
The function $r_{\alpha}^0=r^0(\alpha):\Pi\to \{0,1\}$ of the variable $\alpha=(\alpha_1,\alpha_2)$
is referred to as the {\em initial state} function of the Preisach operator. It
is assumed to be measurable and satisfy the constraints \eqref{v1}, \eqref{v2}, in which case the initial state-input pair 
is called {\em compatible}.
These requirements ensure that the integral in \eqref{pre} is well-defined for each $t\ge 0$ and, furthermore, the output $R_0(t):\mathbb{R}_+\to\mathbb{R}$ of the Preisach model is a continuous function of time.

The function \eqref{relay'} with a fixed $t\ge0$ and varying $\alpha\in \Pi$ is interpreted as the {\em state function} of the Preisach model at the moment $t$ because it describes the states of all the relays at this moment; this state function $r_\alpha(t)=r(t;\alpha):\Pi\to\{0,1\}$ is an element of the space $L_1(\Pi;\mathbb{R})$ for each $t\ge0$.
In this sense, formula \eqref{relay'} defines the evolution of the state function 
in the space $L_1(\Pi;\mathbb{R})$ for a given input provided that the initial state-input pair is compatible.
Moreover, the state function $r_\alpha(t)=r(t;\alpha):\Pi\to\{0,1\}$ and the concurrent input value $I(t)$ are compatible for each $t\ge0$ and the {\em semigroup property} holds:
\begin{equation}\label{semif}
 r_\alpha(\tau)=({\mathcal R}_{\alpha}[r^0_\alpha]I)(\tau) \ \ \ \Rightarrow \ \ \ r_\alpha(t)=({\mathcal R}_{\alpha}[r_\alpha(\tau)]I)(t)\ \ \ \text{for any}\ \ \  \alpha\in\Pi,\ t\geq \tau\geq 0.
\end{equation}

For brevity, let us denote the input-to-output operator of the Preisach model defined by \eqref{pre} as
\begin{equation}\label{P}
R_0(t)=({\mathcal P}[r^0_\alpha]I)(t),\qquad t\geq 0,
\end{equation}
where both the input $I(t):\mathbb{R}_+\to \mathbb{R}$ and the initial state function $r^0_\alpha=r^0({\alpha}): \Pi\to\{0,1\}$ (which is compatible with the initial value of the input)
are the arguments; the value of this operator is the output
$R_0(t):\mathbb{R}_+\to \mathbb{R}$. 
%
%
%
%
An important property of the Preisach operator \eqref{pre} is that it is 
Lipschitz continuous if $q(\alpha):\Pi\to\mathbb{R}_+$ is bounded (Krasnosel'skii et al. 1989). More precisely, the relations 
\[
R_0^k(t)=({\mathcal P}[r^{0,k}_\alpha]I^k)(t),\qquad t\geq 0, {\qquad k=1,2,}
\]
and $0\le I^1(t), I^2(t)\le 1$ ($t\ge 0$) imply
\begin{equation}\label{LipP*}
\|R_0^1(t)-R_0^2(t)\|_{C([0,\tau];\mathbb{R})} \le  q_0 \Big( \|r^{0,1}_\alpha-r^{0,2}_\alpha\|_{L_1(\Pi;\mathbb{R})} + \|I^1(t)-I^2(t)\|_{C([0,\tau];\mathbb{R})}\Big)
\end{equation}
for any $\tau\geq 0$ with the Lipschitz constant
\begin{equation}\label{K}
q_0 := (R_0^{nat}-R_0^{int})\, \sup_{
\alpha\in \Pi} q(\alpha).
\end{equation}

\medskip
\noindent
{\bf 5.2. Phase space of the system with the Preisach operator.}
Let $\mathfrak U$ be the set of all triplets $(I_0,S_0,r_\alpha^0)$, where $(I_0,S_0)\in \mathfrak D=\{(I,S): 0\le I, S; I+S\le 1\}$ 
	and the state function $r^0_\alpha=r^0(\alpha):\Pi\to\{0,1\}$ of the Preisach operator is compatible with $I_0$. 
	The set $\mathfrak U$  is a natural phase space for system \eqref{r0model}, \eqref{pre}.
	In particular, the global Lipschitz estimate \eqref{LipP*} ensures (for example, using the Picard-Lindel\"of type of argument) that for a
	given $(I_0,S_0,r^0_\alpha)\in \mathfrak U$, 
	system \eqref{r0model}, \eqref{pre} has a unique local solution
	with the initial data $I(0)=I_0, S(0)=S_0$ and the initial state function $r^0_\alpha$ of the Preisach opeartor (see e.g. the survey Leonov et al. 2017). This solution satisfies $(I(t),S(t),r_\alpha(t))\in \mathfrak U$ for each $t$ in its domain. Further, the positive invariance of $\mathfrak D$ implies that each solution is extendable to the whole semi-axis $t\in \mathbb{R}_+$. 	Due to the semigroup property \eqref{semif}, these solutions induce a continuous {semi-flow} in the phase space $\mathfrak U$, which is 
	endowed with a metric by the natural embedding into the space $\mathbb{R}^2\times L_1(\Pi;\mathbb{R})$. This construction leads to the standard definition of {\em local and global stability} including stability of equilibrium states and periodic solutions. In particular, an equilibrium is a triplet $(I_0,S_0,r^0_\alpha)\in \mathfrak{U}$ and a periodic solution is a periodic function $(I(t),S(t),r_\alpha(t)): \mathbb{R}_+\to \mathfrak U$ where the last component, viewed as a function $r_\alpha(t)=r(t;\alpha):\mathbb{R}_+\times \Pi\to \{0,1\}$ of two variables $t\in\mathbb{R}_+$ and $\alpha\in\Pi$, is given by \eqref{re}.
The basic reproduction number \eqref{pre} at an equilibrium is constant, while for a periodic solution the basic reproduction number $R_0(t)$ and the state function of the Preisach operator change periodically with the period of $I(t)$ and $S(t)$.

\medskip
\noindent
{\bf 5.3. Switched system. }  A fundamental property of hysteresis operators, which is often used as their definition (see, e.g.~Visisntin 2006; Kluwer Encyclopedia of Mathematics 1997),
 is {\em rate-independence}.  It means that for any increasing continuous function $\theta(t):\mathbb{R}_+\to \mathbb{R}_+$ (transformation of time) inputs $I(t)$ and $I(\theta(t))$ produce identical input-output diagrams such as in Figure \ref{fig_loops}.
Due to the rate-independence property of the Preisach operator,
given an initial state function $r^0=r^0_\alpha:\Pi\to \{0,1\}$, there is a continuous function $ R_{r^0}: \mathbb{R}\to \mathbb{R}$
such that every {\em monotone} input $I(t): \mathbb{R}_+\to\mathbb{R}$ (which is compatible with the state function at the initial moment) and the corresponding output
\eqref{P} of the Preisach operator are related by
\[
R_0(t)= R_{r^0}(I(t)), \qquad t\ge 0.
\]
The function $R_{r^0}(I)$ will be referred to as a {\em branch} of the Preisach operator.
Hence, the state function $r^0=r^0_\alpha$ parameterizes the infinite set of all possible branches.
Due to the assumption that the density $q=q(\alpha)$ of the Preisach operator is positive,
every branch $R_{r^0}(I)$ is a decreasing function on the interval $0\le I\le 1$.
Along with $R_{r^0}(I)$, the functions
\begin{equation}\label{f}
f_{r^0} (I)= I R_{r^0}(I)
\end{equation}
will play an important role.

With the above notation, if the first component $I(t)$ of a solution $(I(t),S(t), r(t)): \mathbb{R}_+\to \mathfrak{U}$ of system \eqref{r0model}, \eqref{pre}
is monotone on a time interval $[t_k,t_{k+1}]$, then the pair $(I(t),S(t))$ is a solution of system
\begin{equation}\label{branch}
\begin{aligned}
    \dot I &=  R(I) SI - I, \\
    \dot S &= -  R(I) SI - \rho S + \rho
\end{aligned}
\end{equation}
with
\begin{equation}\label{branch1}
    R(I)=R_{r^k}(I),
\end{equation}
where $r^k=r(t_k)$, on the same time interval. This ordinary differential system will be called a branch 
of the operator-differential system \eqref{r0model}, \eqref{pre}; again, the state function $r^k=r^k_\alpha$ of the Preisach operator
parameterizes the set of ordinary differential systems \eqref{branch}, \eqref{branch1}.

We see that dynamics of system \eqref{r0model}, \eqref{pre} can be interpreted as switched dynamics of planar systems 
\eqref{branch}, \eqref{branch1}, where the switching moments are the turning points of the component $I$ of a trajectory.
As such, a switch occurs when a trajectory $(I(t),S(t))$ of system \eqref{branch}, \eqref{branch1} (with a fixed $r^k$)
reaches the nullcline $\dot I=0$, i.e.\ at the moment when $ R_{r^k}(I(t)) S(t)=1$ (notice that the nullcline depends on $r^k$).

In what follows, the pair $(I(t),S(t))$ along with the triplet $(I(t),S(t),r(t)): \mathbb{R}_+\to \mathfrak{U}$  will be referred to as a solution of system \eqref{r0model}, \eqref{pre}
whenever this doesn't cause confusion.

\medskip
\noindent
{\bf 5.4. Proofs. }
Throughout the proofs we assume that $R_0^{nat}-R_0^{int}>0$ is sufficiently small
to ensure that the quantities $\epsilon_0, \kappa$ defined below by \eqref{eps0}, \eqref{kappa} using \eqref{K}, \eqref{estim} satisfy 
$\epsilon_0, \kappa>0$.

We first recall well-known stability properties of the planar ordinary differential system \eqref{branch}, 
which are ensured by the following lemmas.

   \begin{lemma}\label{l33} 
   If 
   \begin{equation}\label{eps0}
       \epsilon_0:=R_0^{int}-q_0>0
   \end{equation}
   with $q_0$ given by \eqref{K}, then
   for any branch $R(I)$ of the Preisach operator, the derivative of the function 
   $f(I)=I R(I)$ satisfies
\[
       f'(I)\ge \epsilon_0. 
\]
   \end{lemma}


{\em Proof. } From \eqref{LipP*}, it follows that $-q_0\le R'(I)\le 0$. Combining this with $R(I)\ge R_0^{int}$ and $I\le 1$, we obtain $f'(I)=R(I)+I R'(I)\ge \epsilon_0$. \hfill $\Box$
    
   \begin{lemma}\label{l3}  
   System \eqref{branch}, where $R(I)$ is any branch of the Preisach operator, has a unique infection-free equilibrium $(I^*,S^*) = (0,1)$ 
  and a unique positive (endemic) equilibrium $(I_*,S_*)$, 
  where $I_*$ is a unique solution of
  \begin{equation}\label{nuli}
  \frac{1}{R(I)}=1-\frac{I}{\rho}
  \end{equation}
  and $S_*=1/R(I_*)$ (cf.\ \eqref{infectionfree}, \eqref{endemic}).
  The infection-free equilibrium is unstable, while the positive equilibrium $(I_*,S_*)$ is globally stable in the positive quadrant $I, S>0$.
     \end{lemma}
     
{\em Proof. }
    By assumption, $1/R(0)=1/R_{nat}<1$ and $1/R(\rho)>0$, hence equation \eqref{nuli}
    has a solution $I_*$ in the interval $[0,\rho]\subset[0,1]$. 
    The measure density $q(\alpha)$ of the Preisach operator \eqref{pre} is assumed to be positive,
    which implies that $R(I)$ is a decreasing function.
    Since the function $1/R(I)$ increases, while the function
    $1-I/\rho$ decreases, this solution is unique.
    
    The Jacobi matrix of system \eqref{branch} at the infection-free equilibrium $(0,1)$
    has eigenvalues $\lambda_1=R_{nat}-1, \lambda_2=-R_{nat}-\rho $, hence this equilibrium is a saddle due to the assumption
   $R_{nat}>1.$

    At the endemic equilibrium $(I_*,S_*)$, the characteristic 
    polynomial of the Jacobi matrix of system \eqref{branch} equals
    $$P(\lambda)=\lambda^2+\lambda\bigl(I_* R(I_*)+\rho-I_* R'(I_*)S_*\bigr)+I_* R(I_*)-\rho I_* R'(I_*)S_*.$$
Since $R(I)$ decreases,  all  coefficients of this polynomial are positive,
hence all eigenvalues have negative real parts, and the endemic equilibrium is asymptotically stable.
    
We claim that 
    \begin{equation}\label{V}
    V(I,S)=\int_{I_*}^{I}\left(1-\frac{f(I_*)}{f(i)}\right)\,di+\int_{S_*}^S \left(1-\frac{S_*}{s}\right)\,ds
    \end{equation}
    \[
    = I-I_*+S-S_*-f(I_*)\int_{I_*}^{I}\frac{di}{f(i)} -S_*\ln{\frac{S}{S_*}}.
    \]
is a Lyapunov function for 
system \eqref{branch} in the positive quadrant. 

Indeed, 
 $V(I_*,S_*)=0$, and Lemma \ref{l33} implies that $V(I,S)>0$ for 
 $I, S>0$ except at the point $(I_*,S_*)$. Further, $V\to \infty$ if either $I+S\to \infty$ or $S\to 0$
 or $I\to 0$ (the latter follows from $f(0)=0$).
On the other hand, using the notation 
\[
R_*=R(I_*), \qquad f_*=f(I_*)=I_* R(I_*)
\]
and $R=R(I)$, $f=f(I)$,
the  
 time derivative of $V$ along a trajectory of system \eqref{branch} equals
    \[
\begin{array}{lll}
    \dot V&=& \displaystyle \left(1-\frac{S_*}{S}\right)\rho(1-S)-I+S_*f+\frac{f_* I}{f}-Sf_*\\
&=& \displaystyle \frac{\rho}{S}(S-S_*)(1-S)-f_*(S-S_*)+S_*(f-f_*)-I\left(1-\frac{f_*}{f}\right)\\
    &=& \displaystyle \frac{\rho}{S}(S-S_*)(1-S)-(S-S_*)\frac{\rho(1-S_*)}{S_*}+(S_*f-I)\left(1-\frac{f_*}{f}\right),
       \end{array}
       \]
    where we used that $I_*=\rho (1-S_*)=\rho (1-I_*/f_*)$ due to Lemma \ref{l3}. Equivalently,
     \[
\begin{array}{lll}  
\dot V &=& \displaystyle
    \rho(S-S_*)\left(\frac{1}{S}-\frac{1}{S_*}\right)+\left(\frac{I_*f}{f_*}-I\right)\left(1-\frac{f_*}{f}\right)\\
    &=& \displaystyle -\frac{\rho}{SS_*}(S-S_*)^2+\frac{I_*}{f}(f-f_*)\cdot\frac{I}{f_*}\left(\frac{f}{I}-\frac{f_*}{I_*}\right)\\
    &=& \displaystyle -\frac{\rho}{SS_*}(S-S_*)^2+\frac{1}{RR_*}(f-f_*)(R-R_*). 
   \end{array}
\] 
Since the function $R=R(I)$ decreases and the function $f=f(I)$ increases (Lemma \ref{l33}), we conclude that
\begin{equation}\label{Vadd}
    \dot V\le -\frac{\rho}{SS_*}(S-S_*)^2<0, \qquad (I,S)\ne (I_*,S_*).
\end{equation}
Hence, the equilibrium $(I_*,S_*)$ is globally stable
in the positive quadrant. \hfill $\Box$

\medskip

It turns out that the function $V$ is convex.

\begin{lemma}\label{l1}  
 All the level sets $\{(I,S): V(I,S)<c\}$, $c>0$, of the function $ V$ are convex.
\end{lemma}

{\em Proof. }
Denoting the partial derivatives of the function $V$ by $V_I, V_S$, the curvature of the level line $\{(I,S): V(I,S)=c\}$  of $V$ for any $c>0$ is given by 
\[
\kappa = - \frac{V_{SS} (V_I)^2 + V_{II} (V_S)^2 - 2 V_{IS} V_I V_S}{(V_S^2 + V_I^2)^{3/2}}= -\frac{V_{SS}(V_I)^2+V_{II}(V_S)^2}{(V_S^2+V_I^2)^{3/2}},\]
because $V_{IS}= 0$.
%
Using that
\[
V_I=1-\frac{f(I_*)}{f(I)}, \quad V_S=1-\frac{S_*}{S}, \quad
    V_{II}=\frac{f'(I) f(I_*)}{f^2(I)}, \quad V_{SS}=\frac{S_*}{S^2}
\]
and that $f'(I)>0$ due to Lemma \ref{l33}, we obtain $\kappa<0$, which completes the proof. \hfill $\Box$
    
    
    
  
\medskip
    
Let us consider the trajectory $(I(\cdot),S(\cdot),r(\cdot)): \mathbb{R}_+\to \mathfrak{U}$ of system \eqref{r0model}, \eqref{pre} and the corresponding planar curve
$(I(\cdot),S(\cdot)): \mathbb{R}_+\to \mathbb{R}^2$. Consider the sequence $\mathbb{T}=\{t_k\}$
of successive switching moments $0=t_0<t_1<t_2<\cdots$  such that

\begin{itemize}
    \item 
The restriction of $(I(\cdot),S(\cdot)): \mathbb{R}_+\to \mathbb{R}^2$ to each interval $(t_k,t_{k+1})$, $t_k, t_{k+1}\in \mathbb{T}$ is a solution of system \eqref{branch} with $R(I)=R_{r(t_k)}(I)$
satisfying
\begin{equation}\label{i1}
R_{r(t_k)}(I(t)) S(t)\ne 1, \quad t\in (t_k,t_{k+1}); \qquad R(I(t_{k+1})) S(t_{k+1})= 1.
\end{equation}

 \item
 In addition, if the set $\mathbb{T}=\{t_0,\ldots, t_N\}$ is finite, then the restriction of $(I(\cdot),S(\cdot)): \mathbb{R}_+\to \mathbb{R}^2$ to the infinite interval $(t_N,\infty)$ is a solution of system \eqref{branch} with $R(I)=R_{r(t_N)}(I)$
satisfying
\[
R_{r(t_N)}(I(t)) S(t)\ne 1, \quad t>t_N.
\]
\end{itemize}

If the set $\mathbb{T}$ is finite, the  trajectory of system \eqref{r0model}, \eqref{pre} follows the flow \eqref{branch}, \eqref{branch1} on the interval $[t_N,\infty)$ and due to Lemma \ref{l3} converges to an equilibrium. In this scenario, this equilibrium of \eqref{branch}, \eqref{branch1} is a node.

We continue the proof assuming that $\mathbb{T}=\{t_k\}$ is infinite, i.e.\ there is an infinite sequence of switches along the trajectory.
Let 
\[
I_k=I(t_k), \quad S_k=S(t_k).
\]
Denote by $V_k(I,S)$ the Lyapunov function \eqref{V} for system \eqref{branch}, \eqref{branch1} (where $r^k=r(t_k)$) and by $(I_*^k,S_*^k)$ the positive (endemic) equilibrium of this system. We will also use the notation $R_k(I)=R_{r(t_k)}(I)$.
The next lemma provides an estimate for the (negative) increment 
\begin{equation}\label{deltaV}
\Delta V_k :=V_k(I_{k+1},S_{k+1})-V_k(I_k,S_k)
\end{equation}
of this function along the trajectory on the time interval $[t_k,t_{k+1}]$.

\begin{figure}[H]
\centering 
\includegraphics[width=1.0\textwidth]{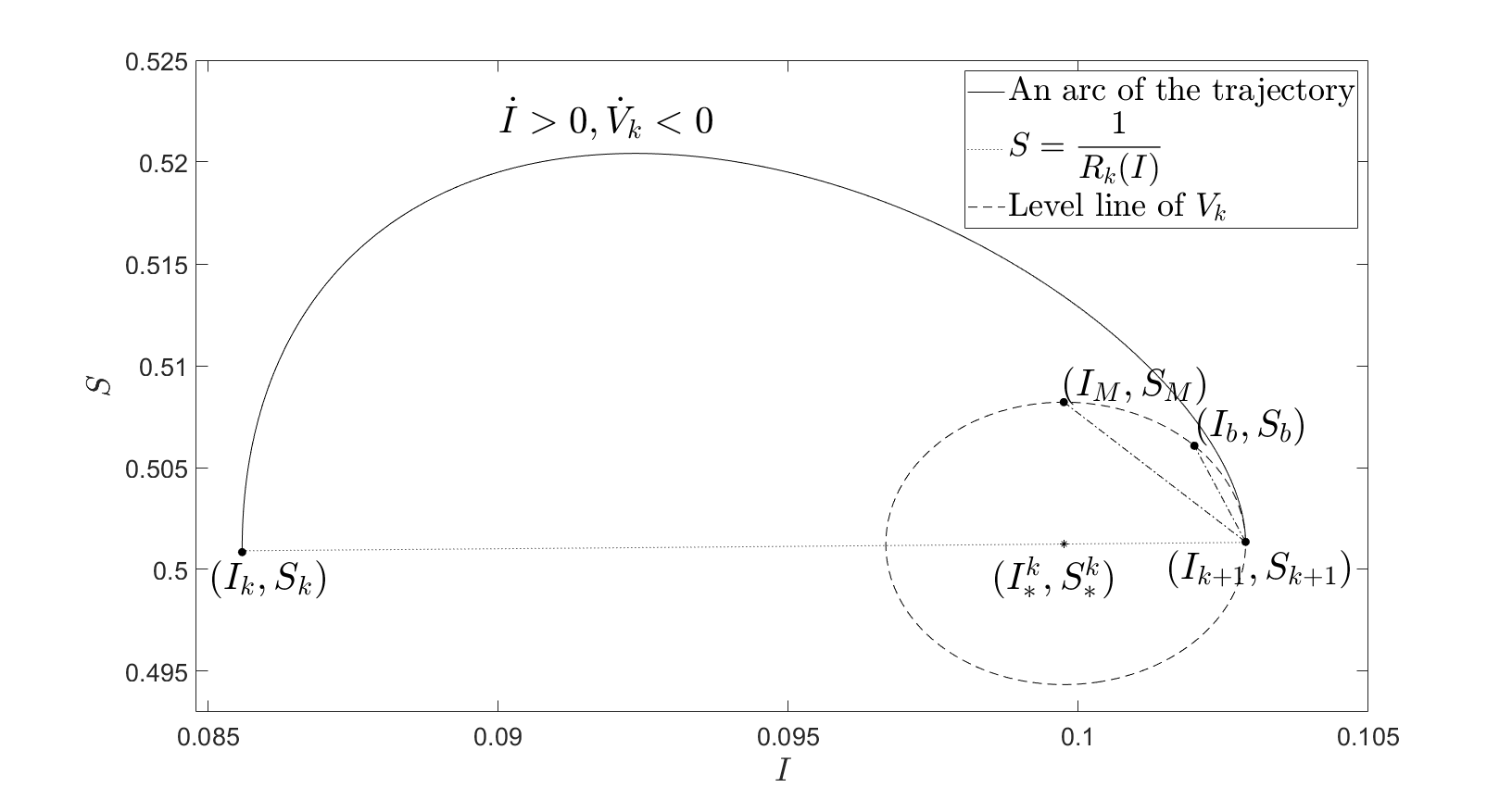} 
\caption{An arc of a trajectory of system \eqref{r0model}, \eqref{pre} satisfying $R_k(I(t)) S(t)\ge 1$. \label{fig1}}
\end{figure}

\begin{lemma}\label{l2}  
If $R_{r(t_k)}(I(t)) S(t)> 1$ on $(t_k,t_{k+1})$ (cf.\ \eqref{i1}), then the increment \eqref{deltaV} of the Lyapunov function $V_k$ along the trajectory 
satisfies
    \begin{equation}\label{VV}
    \Delta V_k\le -\frac{\rho}{4}(I_{k+1}-I_*^k)\sqrt{S_M\int_{I_*^k}^{I_{k+1}}\left(1-\frac{f(I_*^k)}{f(i)}\right)di},
    \end{equation}
    where     $M=(I_M,S_M)$ is the highest point of the level line $V_k(I,S)=V_k(I_{k+1},S_{k+1})$
of the Lyapunov function $V_k$, see Figure \ref{fig1}.
\end{lemma}    

{\em Proof.} 
    At the point $M=(I_M,S_M)$,
    $$0=\left.\frac{\partial V_k}{\partial I}\right|_M=1-\frac{f(I_*^k)}{f(I_M)},$$ which is equivalent to $f(I_*^k)=f(I_M)$. Since $f(I)$ increases, it follows that  $I_M=I_*^k$,
    hence
     $M=(I_*^k,S_{M})$ and
    $$V_k(I_*^k,S_{M})=S_{M}-S_*^k-S_*^k\ln{\frac{S_{M}}{S_*^k}}.$$
    
    Consider a point  $(I_b,S_b)$ on the arc of the level line of $V_k$ connecting the points $M$ and 
    $(I_{k+1},S_{k+1})$.
    Due to Lemma \ref{l1},
    $$\frac{S_b-S_{k+1}}{I_{k+1}-I_b}\ge \frac{S_{M}-S_{k+1}}{I_{k+1}-I_*^k}$$
    (see Figure \ref{fig1}). Equivalently,
    $$I_{k+1}-I_b\le \frac{S_b-S_{k+1}}{S_{M}-S_{k+1}}(I_{k+1}-I_*^k).$$
    This implies the estimate
    \begin{equation}\label{bbb0}
         I_b-I_*^k\ge I_{k+1}-I_*^k-\frac{S_b-S_{k+1}}{S_{M}-S_{k+1}}(I_{k+1}-I_*^k)=\frac{S_{M}-S_b}{S_{M}-S_{k+1}}(I_{k+1}-I_*^k).
    \end{equation}
    In the domain $I\ge I_*^k, S\ge {1}/{R_k(I)},$ we have
    ${1}/{R_k(I)}\ge {1}/{R_k(I_*^k)}=S_*^k,$ which implies 
    $0\le S-{1}/{R_k(I)}\le S-S_*^k$.
    Hence, relation \eqref{Vadd} implies that along the trajectory,
    \begin{eqnarray}
    \nonumber
    \frac{dV_k}{dI}=\frac{\dot V_k}{\dot I}\le -\frac{\rho}{SS_*^kIR_k(I)}\cdot\frac{(S-S_*^k)^2}{(S-1/R_k(I))}
    \le 
    -\frac{\rho}{SS_*^k I R_k(I_*^k)}(S-S_*^k)
    \\
    \nonumber
    =-\frac{\rho}{SI}(S-S_*^k)\le -\frac{\rho}{I}(S-S_*^k).
    \end{eqnarray}
    Since $V_k$ decreases along the trajectory,
    the segment of the trajectory with $t\in(t_k,t_{k+1})$ lies above the level line of $V_k$, hence for this segment 
    $I\in [I_*^k,I_b]$ implies
    $S\ge S_b$  
and
    \begin{equation}\label{vid}
        \frac{dV_k}{dI}\le -\frac{\rho}{I}(S_b-S_*^k).
    \end{equation}
    The relation $\dot I(t) > 0$, $t\in(t_k,t_{k+1})$, implies that the arc 
    of the trajectory for this time interval admits a parameterization 
    $(I,\hat S(I))$, $I\in [I_k,I_{k+1}]$.
    Therefore, integrating \eqref{vid} over the interval $[I_k^*,I_b]$, we obtain
    $$V_k(I_b,\hat S(I_b))- V_k(I_*^k,\hat S(I_*^k)) \le -\rho(S_b-S_*^k)\ln{\frac{I_b}{I_*^k}}.$$
    Further, using the identity $\ln (1+x)\ge {x}/(1+x),$ 
    $x\ge -1$, we arrive at
    \[
    V_k(I_b,\hat S(I_b))- V_k(I_*^k,\hat S(I_*^k))  
    \le 
    -\rho(S_b-S_*^k)\frac{I_b-I_*^k}{I_b}\le -\rho(S_b-S_*^k)(I_b-I_*^k)
 \]
and due to \eqref{bbb0},
    \begin{equation}\label{V}
V_k(I_b,\hat S(I_b))- V_k(I_*^k,\hat S(I_*^k))  \le -\rho\frac{(S_b-S_*^k)(S_{M}-S_b)(I_{k+1}-I_*^k)}{S_{M}-S_{k+1}}.
    \end{equation}
    
The maximum of the function  $\frac{(S_b-S_*^k)(S_{M}-S_b)}{S_{M}-S_{k+1}}$ on the interval $[S_{k+1},S_{M}]$ is achieved at $S_{k+1}= (S_{M}+S_*^k)/{2},$ 
hence
$$
\max_{S_b\in [S_{k+1},S_{M}]}\frac{(S_b-S_*^k)(S_{M}-S_b)}{S_{M}-S_{k+1}}\geq
\frac{(S_{M}-S_*^k)^2}{4(S_{M}-S_{k+1})}
\ge \frac{(S_{M}-S_*^k)^2}{4(S_{M}-S_*^k)}
=\frac{S_{M}-S_*^k}{4}.
$$
    
    Therefore, from \eqref{V} it follows that 
    $$ V_k(I_b,\hat S(I_b))- V_k(I_*^k,\hat S(I_*^k))  \le -\frac{\rho}{4}(I_{k+1}-I_*^k)(S_{M}-S_*^k).$$ 
    Since $V_k$ decreases along the trajectory, this estimate implies 
       \begin{equation}\label{l4*}
    \Delta V_k\le-\frac{\rho}{4}(I_{k+1}-I_*^k)(S_{M}-S_*^k).
    \end{equation}

Finally, let us estimate $S_M-S_*^k.$ Since $V_k$ has the same value 
at the points $(I_*^k, S_M)$ and $(I_{k+1}, S_{k+1})$, we have 
    $$S_{M}-S_*^k-S_*^k\ln{\frac{S_{M}}{S_*^k}}=I_{k+1}-I_*^k-f(I_*^k)\int_{I_*^k}^{I_{k+1}}\frac{di}{f(i) }+S_{k+1}-S_*^k-S_*^k\ln{\frac{S_{k+1}}{S_*^k}},$$\\
    which implies
    \[
    S_{M}-S_{k+1}-S_*^k\ln{\frac{S_{M}}{S_{k+1}}}=\int_{I_*^k}^{I_{k+1}}\left(1-\frac{f(I_*^k)}{f(i)}\right)di.
    \]
    We estimate the left hand side of this equation using $\ln{(1+x)}\ge {x}/{1+x}$, $x\ge-1$:
    \[
    S_{M}-S_{k+1}-S_*^k\ln{\frac{S_{M}}{S_{k+1}}}
    \le 
    S_{M}-S_{k+1}-\frac{S_*^k(S_{M}-S_{k+1})}{S_{M}}=\left(1-\frac{S_*^k}{S_{M}}\right)\left(S_{M}-S_{k+1}\right),
    \]
    hence
    $$
    \frac{1}{S_{M}}(S_{M}-S_*^k)(S_{M}-S_{k+1})\ge \int_{I_*^k}^{I_{k+1}}\left(1-\frac{f(I_*^k)}{f(i)}\right)di. 
    $$
 Since $S_*^k< S_{k+1}$, 
  this implies
     \begin{equation}\label{SS}
    (S_{M}-S_*^k)^2\ge S_{M}\int_{I_*^k}^{I_{k+1}}\left(1-\frac{f(I_*^k)}{f(i)}\right)di.
    \end{equation}
Combining this relation with \eqref{l4*}, we obtain \eqref{VV}.\hfill $\Box$

\medskip
Next, we estimate the difference
\begin{equation}\label{tildeV}
\tilde \Delta V_k :=V_{k+1}(I_{k+1},S_{k+1})-V_k(I_{k+1},S_{k+1})
\end{equation}
between the values of  Lyapunov functions $V_{k+1}$ and $V_k$ at the point $(I_{k+1},S_{k+1})$. We begin with the following statement.
    

\begin{lemma}\label{l5'} 
Under the assumptions of Lemma \ref{l2}, 
\begin{equation}\label{lemt}
0\le I_*^k-I_*^{k+1} = \rho(S_*^{k+1}-S_*^k)\le \frac{\rho q_0}{R_{k+1}(I_*^{k+1})R_k(I_*^k)}\,(I_{k+1}-I_*^k ).
\end{equation}
\end{lemma} 
   
{\em Proof. } 
Consider the branch $R_k(I):=R_{r(t_k)}(I)$  of the Preisach operator (see Section 5.3), the corresponding function $f_k(I):= I R_k(I)$ 
and the endemic equilibrium $(I_*^k,S_*^k)$ of system \eqref{branch}, \eqref{branch1} for this branch:
 $$
 S_*^k=1-\frac{I_*^k}{\rho}=\frac{1}{R_k(I_*^k)}.
 $$
These equations imply
    $$ S_*^{k+1}-S_*^k=\frac{1}{\rho}(I_*^k-I_*^{k+1})$$
    and 
    \begin{eqnarray}\nonumber
    \frac{I_*^k-I_*^{k+1}}{\rho}=\frac{1}{R_{k+1}(I_*^{k+1})}-\frac{1}{R_k(I_*^k)}
    =\frac{R_k(I_*^k)-R_{k+1}(I_*^k)+R_{k+1}(I_*^k)-R_{k+1}(I_*^{k+1})}{R_{k+1}(I_*^{k+1})R_k(I_*^k)}\\\nonumber
   =\frac{1}{R_{k+1}(I_*^{k+1})R_k(I_*^k)}\bigr(R_k(I_*^k)-R_{k+1}(I_*^{k}) + R_{k+1}'(\hat I)(I_*^k-I_*^{k+1})\bigr)\\\nonumber
    \end{eqnarray}
    for an intermediate value of $\hat I$ between $I_*^k$ and $I_*^{k+1}$.
    
    If we assume that $I_*^{k+1}>I_*^k$, then we arrive at a contradiction:
    $$0>\frac{1}{\rho}(I_*^k-I_*^{k+1})=\frac{1}{R_{k+1}(I_*^{k+1})R_k(I_*^k)}\bigr(R_k(I_*^k)-R_{k+1}(I_*^{k}) + R_{k+1}'(\hat I)(I_*^k-I_*^{k+1})\bigr)\ge 0$$\\
    because $R_{k+1}(I)$ decreases and \eqref{LipP*} implies that 
    \begin{equation}\label{ddt}
   0\le  R_k(I_*^k)- R_{k+1}(I_*^{k}) \le q_0( I_{k+1}-I_*^k).
    \end{equation}
    Hence, $I_*^{k+1}\le I_*^k$, and consequently  $S_*^{k+1}\ge S_*^k$
   as well as
    \begin{equation}\label{RR}
  0\le S_*^{k+1}-S_*^k= \frac1\rho( I_*^k-I_*^{k+1}) \le \frac{R_k(I_*^k)-R_{k+1}(I_*^{k})}{R_{k+1}(I_*^{k+1})R_k(I_*^k)}.
        \end{equation}
%
 Combining \eqref{ddt} and \eqref{RR}, we obtain \eqref{lemt}.\hfill $\Box$
 
 \medskip
 We use Lemma \ref{l5'} to prove the following estimate for the difference \eqref{tildeV}.
 
 \begin{lemma}\label{l5} 
Under the assumptions of Lemma \ref{l2}, the difference \eqref{tildeV} satisfies
\begin{equation}\label{eql5}
\tilde \Delta V_k \le 
\frac{q_0  (I_{k+1}-I_*^k)^2}{R_k(I_*^k)}  \left(\frac{\rho^2 q_0}{R_{k+1}^2(I_*^{k+1})f_k(I_*^k) }+\frac{\rho}{f_k(I_*^k)}+1\right).
\end{equation}
\end{lemma} 
    
 {\em Proof.}   
    By the definition of the Lyapunov function, $\tilde \Delta V_k=\Delta _S+ \Delta _I$, where
    \begin{eqnarray}\nonumber
  \Delta _S&=&  \int_{S_*^{k+1}}^{S_{k+1}}\left(1-\frac{S_*^{k+1}}{s}\right)\,ds-\int_{S_*^k}^{S_{k+1}}\left(1-\frac{S_*^k}{s}\right)\,ds
  \\ \nonumber
    &=&  -\int_{S_*^k}^{S_*^{k+1}}\left(1-\frac{S_*^k}{s}\right)\,ds-\int_{S_*^{k+1}}^{S_{k+1}}\left(\frac{S_*^{k+1}-S_*^k}{s}\right)\,ds,
     \end{eqnarray}
         \begin{eqnarray}\nonumber
\Delta_I&=&    \int_{I_*^{k+1}}^{I_{k+1}}\left(1-\frac{f_{k+1}(I_*^{k+1})}{f_{k+1}(i)}\right)\,di-\int_{I_*^k}^{I_{k+1}}\left(1-\frac{f_k(I_*^k)}{f_k(i)}\right)\,di \\\nonumber
    &=&-\int_{I_*^k}^{I_{k+1}}\left(\frac{f_{k+1}(I_*^{k+1})}{f_{k+1}(i)} 
    -\frac{f_k(I_*^k)}{f_k(i)}\right)\,di
    -\int_{I_*^k}^{I_*^{k+1}}\left(1-\frac{f_{k+1}(I_*^{k+1})}{f_{k+1}(i)}\right)\,di \\ \nonumber &=&-
    \int_{I_*^k}^{I_{k+1}}\left(\frac{f_{k+1}(I_*^{k+1})}{f_{k+1}(i)} -\frac{f_{k+1}(I_*^{k+1})}{f_k(i)}\right)\,di
   -\int_{I_*^k}^{I_{k+1}}\frac{f_{k+1}(I_*^{k+1})-f_k(I_*^k)}{f_k(i)}\,di
   \\ \nonumber
    &-&\int_{I_*^k}^{I_*^{k+1}}\left(1-\frac{f_{k+1}(I_*^{k+1})}{f_{k+1}(i)}\right)\,di.
    \end{eqnarray}
By Lemma \ref{l5'}, $S_*^k\le S_*^{k+1}\le S_{k+1}$, hence 
 \begin{equation} \label{Spart}
 \int_{S_*^k}^{S_*^{k+1}}\left(1-\frac{S_*^k}{s}\right)\,ds\ge 0, \quad \int_{S_*^{k+1}}^{S_{k+1}}\left(\frac{S_*^{k+1}-S_*^k}{s}\right)\,ds\ge 0,
\end{equation}
and consequently $\Delta_S\le 0$.
 Also, the relation $R_{k+1}(I)\le R_k(I)$, $I\le I_{k+1}$, implies that
 \begin{equation}\label{kl}
 \int_{I_*^k}^{I_{k+1}}\left(\frac{f_{k+1}(I_*^{k+1})}{f_{k+1}(i)} -\frac{f_{k+1}(I_*^{k+1})}{f_k(i)}\right)\,di=
 \int_{I_*^k}^{I_{k+1}}\frac{f_{k+1}(I_*^{k+1})}{i}\left(\frac{1}{R_{k+1}(i)} -\frac{1}{R_k(i)}\right)\,di\ge 0,
 \end{equation}
    therefore  
      \begin{equation}\label{neest'}
    \tilde \Delta V_k\le \Delta_I\le 
         -\int_{I_*^k}^{I_*^{k+1}}\left(1-\frac{f_{k+1}(I_*^{k+1})}{f_{k+1}(i)}\right)\,di
      -\int_{I_*^k}^{I_{k+1}}\frac{f_{k+1}(I_*^{k+1})-f_k(I_*^k)}{f_k(i)}\,di.
    \end{equation}
 In order to  estimate the first integral in this expression, we recall that $f_{k+1}(I)$ increases, $R_{k+1}(I)$ decreases and $I_*^{k+1} \le I_*^k$ due to Lemma \ref{l5'}, hence
    \begin{eqnarray}\nonumber
    0&\le& -\int_{I_*^k}^{I_*^{k+1}}\left(1-\frac{f_{k+1}(I_*^{k+1})}{f_{k+1}(i)}\right)\,di\le (I_*^k-I_*^{k+1})\left(1-\frac{f_{k+1}(I_*^{k+1})}{f_{k+1}(I_*^k)}\right)\\\nonumber
    &=&(I_*^k-I_*^{k+1})\frac{R_{k+1}(I_*^k)I_*^k-R_{k+1}(I_*^{k+1})I_*^{k+1}}{R_{k+1}(I_*^k)I_*^k}\le (I_*^k-I_*^{k+1})\frac{R_{k+1}(I_*^k)(I_*^k-I_*^{k+1})}{R_{k+1}(I_*^k)I_*^k}\\ \label{I=I}
    &=&\frac{(I_*^k-I_*^{k+1})^2}{I_*^k}
     \end{eqnarray}
     and using Lemma \ref{l5'} again,
   \begin{equation}\label{neest}  
     -\int_{I_*^k}^{I_*^{k+1}}\left(1-\frac{f_{k+1}(I_*^{k+1})}{f_{k+1}(i)}\right)\,di\le \frac{\rho^2 q_0^2}{R_{k+1}^2(I_*^{k+1})R_k^2(I_*^k) I_*^k}\,(I_{k+1}-I_*^k )^2.
  \end{equation}   
 On the other hand,
    \begin{eqnarray}\nonumber
     &&-\int_{I_*^k}^{I_{k+1}}\frac{f_{k+1}(I_*^{k+1})-f_k(I_*^k)}{f_k(i)}\,di\\
     \nonumber
   && =-\int_{I_*^k}^{I_{k+1}}\frac{f_{k+1}(I_*^{k+1})-f_{k+1}(I_*^k)+f_{k+1}(I_*^k)-f_k(I_*^k)}{f_k(i)}di\\ \nonumber
    &&=-\int_{I_*^k}^{I_{k+1}}\frac{I_*^{k+1}R_{k+1}(I_*^{k+1})-I_*^kR_{k+1}(I_*^k)+I_*^k(R_{k+1}(I_*^k)-R_k(I_*^k))}{f_k(i)}di \\\nonumber
    &&=\int_{I_*^k}^{I_{k+1}}\frac{I_*^{k+1}(R_{k+1}(I_*^{k})-R_{k+1}(I_*^{k+1}))}{f_k(i)}di\\\nonumber
    &&+\int_{I_*^k}^{I_{k+1}}\frac{(I_*^{k}-I_*^{k+1})R_{k+1}(I_*^k)+I_*^k(R_{k}(I_*^k)-R_{k+1}(I_*^k))}{f_k(i)}di\\\label{cou}
    && \le  \frac{I_{k+1}-I_*^k}{f_k(I_*^k)}\bigl((I_*^{k}-I_*^{k+1})R_{k+1}(I_*^k) +I_*^k(R_{k}(I_*^k)-R_{k+1}(I_*^k))\bigr),
  \end{eqnarray}  
  where we used that $R_{k+1}(I_*^{k})-R_{k+1}(I_*^{k+1})\le 0$ because $R_{k+1}(I)$ decreases
  and $f_k(I_*^k)\le f_k(I)$ on $[I_*^k,I_{k+1}]$ because $f_k(I)$ increases.
  Hence, relations  \eqref{lemt} and \eqref{ddt} imply
  \begin{eqnarray}\nonumber
-\int_{I_*^k}^{I_{k+1}}\frac{f_{k+1}(I_*^{k+1})-f_k(I_*^k)}{f_k(i)}\,di\le   
\frac{\rho q_0 R_{k+1}(I_*^k)(I_{k+1}-I_*^k )^2}{R_{k+1}(I_*^{k+1})R_k(I_*^k)f_k(I_*^k)}\\
\nonumber
+ \frac{q_0 I_*^k (I_{k+1}-I_*^k)^2}{f_k(I_*^k)} \le \frac{q_0  (I_{k+1}-I_*^k)^2}{R_k(I_*^k)} \left(\frac{\rho}{f_k(I_*^k)}+1\right).
  \end{eqnarray}
  Combining this estimate with \eqref{neest'} and \eqref{neest}, we obtain \eqref{eql5}. \hfill $\Box$
 
 
\medskip
Now we combine Lemmas \ref{l2} and \ref{l5} to obtain the following statement.
 
 \begin{lemma}\label{lmain}  
If $R_{r(t_k)}(I(t)) S(t)> 1$ on $(t_k,t_{k+1})$, then 
    \begin{equation}\label{VVV'}
    V_{k+1}(I_{k+1},S_{k+1})-V_k(I_k,S_k)
    \le -Q_k(I_{k+1}-I_*^k)^2
    \end{equation}
    with
    \begin{equation}\label{V4'}
    Q_k=\frac{\rho}{4}\sqrt{\frac{\epsilon_0 S_M}{2f_k(I_{k+1})}}
    -\frac{q_0  }{R_k(I_*^k)}  \left(\frac{\rho^2 q_0}{R_{k+1}^2(I_*^{k+1})f_k(I_*^k) }+\frac{\rho}{f_k(I_*^k)}+1\right),
    \end{equation}
    where     $M=(I_M,S_M)$ is the highest point of the level line $V_k(I,S)=V_k(I_{k+1},S_{k+1})$
of the Lyapunov function $V_k$.
\end{lemma}  

{\em Proof.} From Lemma \ref{l33}, it follows that
    \begin{eqnarray}
    \int_{I_*^k}^{I_{k+1}}\left(1-\frac{f_k(I_*^k)}{f_k(i)}\right)di= \int_{I_*^k}^{I_{k+1}}\frac{f_k(i)-f_k(I_*^k)}{f_k(i)}di\\\nonumber\ge 
    \int_{I_*^k}^{I_{k+1}}\frac{\epsilon_0(i-I_*^k)}{f_k(I_{k+1})}di\ge \frac{\epsilon_0(I_{k+1}-I_*^k)^2}{2f_k(I_{k+1})}.
    \end{eqnarray}
    Substituting this estimate in \eqref{VV} and adding \eqref{eql5}, we obtain \eqref{VVV'}. \hfill $\Box$
    

    
    \medskip
As the next step, we obtain a counterpart of Lemma \ref{lmain} for the parts of the trajectory satisfying 
$R_{r(t_k)}(I(t)) S(t) < 1$ on $(t_k,t_{k+1})$. 

 \begin{lemma}\label{lmain'}  
If $R_{r(t_k)}(I(t)) S(t)< 1$ on $(t_k,t_{k+1})$, then 
    \begin{equation}\label{VVV}
    V_{k+1}(I_{k+1},S_{k+1})-V_k(I_k,S_k)
    \le -\tilde{Q_k}(I_{k+1}-I_*^k)^2
    \end{equation}
    with
    \begin{equation}\label{V4}
    \tilde{Q_k}=\frac{\rho}{4}\sqrt{\frac{\epsilon_0 S_m}{2f_k({I_*^k})}}
    -
    q_0\left(
  \frac{\rho^2 q_0}{R_{k+1}^2(I_*^{k+1})R_k^2(I_*^k) I_*^k}+ 
 \frac{\rho R_{k+1}(I_*^k)}{R_{k+1}(I_*^{k+1})R_k(I_*^k)f_k(I_{k+1})}+\frac{I_*^k}{f_k(I_{k+1})}\right),
    \end{equation}
    where     $m=(I_m,S_m)$ is the lowest point of the level line $V_k(I,S)=V_k(I_{k+1},S_{k+1})$
of the Lyapunov function $V_k$.
\end{lemma} 
    
    \begin{figure}[H]
    \centering 
    \includegraphics[width=1\textwidth]{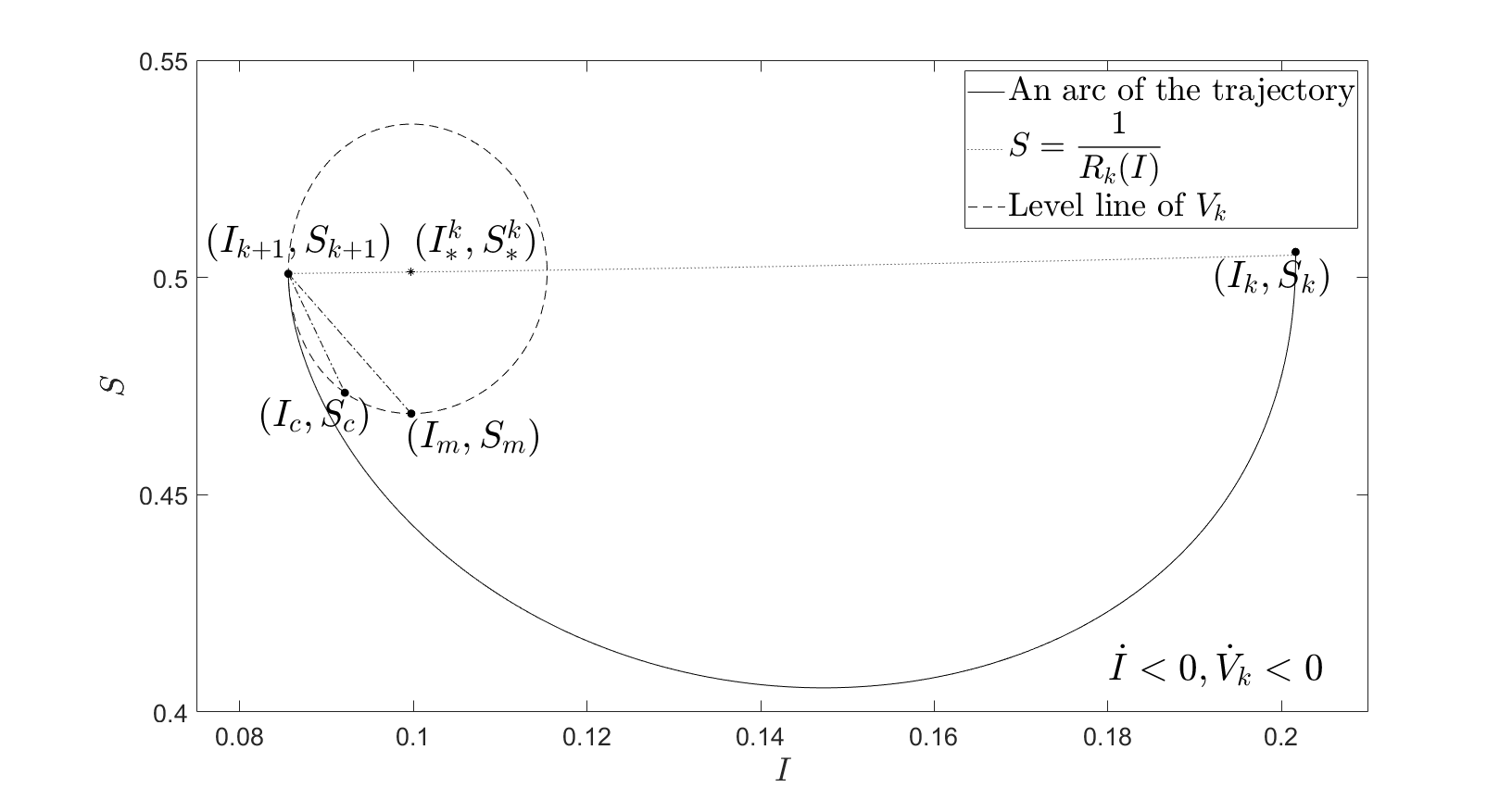} 
    \caption{An arc of a trajectory of system \eqref{r0model}, \eqref{pre} satisfying $R_k(I(t)) S(t)\le 1$. \label{fig2}}
    \end{figure}
    
    {\em Proof. }
We adapt the proofs of Lemmas \ref{l2}\,--\,\ref{l5} 
to obtain their counterparts for the part of the trajectory satisfying $R_{r(t_k)}(I(t)) S(t)< 1$ on $(t_k,t_{k+1})$, which implies $\dot I<0$.
The main steps are outlined below for completeness.
   
The first coordinate of the point $m=(I_m,S_m)$ and the value of the Lyapunov function $V_k$ at this point satisfy
$I_m=I_k^*$ and
\[
    V_k(I_*^k,S_{m})=S_{m}-S_*^k-S_*^k\ln{\frac{S_{m}}{S_*^k}}=V_k(I_{k+1},S_{k+1}),
 \]
 see Figure \ref{fig2}.
 Using Lemma \ref{l1}, we obtain the estimate
    \begin{equation}\label{+}
    I_*^k-I_c\ge I_*^k-I_{k+1}-\frac{S_{k+1}-S_c}{S_{k+1}-S_{m}}(I_*^k-I_{k+1})=\frac{S_c-S_{m}}{S_{k+1}-S_{m}}(I_*^k-I_{k+1})
    \end{equation}
 for any point $(I_c,S_c)$ that belongs to the arc of the level line of $V_k$ connecting the points $m$ and $(I_{k+1},S_{k+1})=(I_{k+1},1/{R_k(I_{k+1})})$
    (cf.\ \eqref{bbb0}).
    
    Using \eqref{Vadd} in the domain $I\le I_*^k$, $S\le {1}/{R_k(I)}$, we obtain that along the trajectory
  \[
    \frac{dV_k}{dI}\ge -\frac{{\rho}(S-S_*^k)^2}{{SS_*^k}(R_k(I)S-1)I}\ge 0,
 \]   
 where
    ${1}/{R_k(I)}\le {1}/{R_k(I_*^k)}=S_*^k $ implies $0\le {1}/{R_k(I)}-S\le S_*^k-S$, hence
    $$
    \frac{dV_k}{dI}\ge \frac{{\rho}(S-S_*^k)^2}{{SS_*^k}R_k(I)I(S_*^k-S)}=\frac{\rho}{SS_*^kR_k(I)I}(S_*^k-S)$$
    and due to $R_k(I)S<1$, $S_*^k\le 1$, 
    $$
    \frac{dV_k}{dI}\ge \frac{\rho(S_*^k-S)}{I}.
    $$
    Since the trajectory lies below the level line of $V_k$, on the interval $[I_c,I_*]$ we have $S\le S_c$ and $\frac{dV_k}{dI}\ge \frac{\rho}{I}(S_*^k-S_c)$,
    which implies 
    $$
    V_k(I_*^k,\hat S(I_*^k))-V_k(I_c, \hat S(I_c))\ge \rho(S_*^k-S_c)\ln{\frac{I_*^k}{I_c}},
    $$
    where $(I,\hat S(I))$, $I\in [I_{k+1},I_k]$ is a parameterization of the segment of the trajectory, and
    using the identity $\ln{(1+x)}\ge x/(1+x)$, $x\ge -1$, we obtain
    \begin{eqnarray}\nonumber
 V_k(I_c, \hat S(I_c)) -  V_k(I_*^k,\hat S(I_*^k))\le 
  -\rho(S_*^k-S_c)\frac{I_*^k-I_c}{I_*^k}\\\nonumber
    \le -\rho\frac{(S_*^k-S_c)(S_c-S_{m})(I_*^k-I_{k+1})}{S_{k+1}-S_{m}},
    \end{eqnarray}
    where in the last step we used $I_*^k\le 1$ and  \eqref{+}. 
 %
 %
 Maximizing the right hand side of this formula over
 $S_c\in [S_m,S_{k+1}]$ as in the proof of Lemma \ref{l2} and taking into account that the function $V_k$ decreases along the trajectory, we arrive at
 (cf.~\eqref{l4*})
    \begin{equation}\label{V3}
  \Delta V_k=  V_k(I_{k+1},S_{k+1})-V_k(I_k,S_k)\le -\frac{\rho}{4}(S_*^k-S_m)(I_*^k-I_{k+1}).
    \end{equation}

    In order to estimate the difference $S_*^k-S_m$ in \eqref{V3}, we 
    equate the values of the function $V_k$ at the points 
    $(I_*^k,S_m)$ and $(I_{k+1},S_{k+1})$ to obtain 
    $$S_*^k\ln{\frac{S_{k+1}}{S_m}}-(S_{k+1}-S_m)=\int_{I_{k+1}}^{I_*^k}\left(\frac{f_k(I_*^k)}{f_k(i)}-1\right)di.$$
    Using the inequality $2(1+x)\ln{(1+x)}\le {x(2+x)}$, $x\ge 0$ we get
    $$\ln{\frac{S_{k+1}}{S_m}}\le 
 \frac{S_{k+1}^2-S_m^2}{2S_{k+1}S_m},$$ 
    where $S_m<S_{k+1}<S_*^k$, therefore
    \begin{eqnarray}\nonumber
 &&   S_*^k\ln{\frac{S_{k+1}}{S_m}}-(S_{k+1}-S_m)\le 
    (S_{k+1}-S_m)\left(\frac{S_*^k(S_{k+1}+S_m)}{2S_{k+1}S_m}-1 \right)
    \\\nonumber
    &&
    =\frac{S_{k+1}-S_m}{2}\left( \frac{S_*^k-S_m}{S_m}+\frac{S_*^k-S_{k+1}}{S_{k+1}}\right)
  \le 
 \frac{(S_*^k-S_m)^2}{S_m},
    \end{eqnarray}
  hence
 $$
 \frac{(S_*^k-S_m)^2}{S_m}\ge \int_{I_{k+1}}^{I_*^k}\left(\frac{f_k(I_*^k)}{f_k(i)}-1\right)di
 $$ 
 and from \eqref{V3} it follows that (cf.~\eqref{VV})
    \[
 \Delta V_k \le -\frac{\rho}{4}(I_*^k-I_{k+1})\sqrt{S_m\int_{I_{k+1}}^{I_*^k}\left(\frac{f_k(I_*^k)}{f_k(i)}-1\right)di}.
    \]
    Further, due to Lemma \ref{l33},
        \[
    \int_{I_*^k}^{I_{k+1}}\left(1-\frac{f_k(I_*^k)}{f_k(i)}\right)di
    \ge \frac{\epsilon_0(I_{k+1}-I_*^k)^2}{2f_k(I_*^{k})},
    \]
    consequently,
    \begin{equation}\label{nes}
   \Delta V_k \le -\frac{\rho}{4}(I_*^k-I_{k+1})^2\sqrt{\frac{\epsilon_0 S_m}{2f_k(I_*^k)}}.      
    \end{equation}

   

As in the proof of Lemma \ref{l5'},
    $$
    S_*^k-S_*^{k+1}=\frac{I_*^{k+1}-I_*^k}{\rho}
    =\frac{1}{R_k(I_*^k)R_{k+1}(I_*^{k+1})}\bigl(R_{k+1}'(\hat I)(I_*^{k+1}-I_*^k)+R_{k+1}(I_*^k)-R_k(I_*^k) \bigr)
    $$
but this time (cf.\ \eqref{ddt})
    \begin{equation}\label{ddt'}
   0\le  R_k(I_*^{k+1})- R_{k}(I_*^{k}) \le q_0( I_*^k-I_{k+1}),
    \end{equation}
which allows us to conclude (due to $R_{k+1}'(\hat I)\le 0 $) that
$
    I_*^{k+1}\ge I_*^k$ and   $S_*^{k+1}\le S_*^k$, 
and the argument used in the proof of Lemma \ref{l5'} leads to the estimates (cf.~\eqref{lemt})
\begin{equation}\label{lemt'}
0\le I_*^{k+1}-I_*^{k} = \rho(S_*^{k}-S_*^{k+1})\le \frac{\rho q_0}{R_{k+1}(I_*^{k+1})R_k(I_*^k)}\,(I_*^k-I_{k+1}).
\end{equation}
In particular, as opposed to the case considered in Lemmas \ref{l2}\,--\,\ref{l5},
\begin{equation}\label{case2}
    I_*^{k+1}\ge I_*^k\ge I_{k+1}, \quad S_{k+1}\le S_*^{k+1}\le S_*^k,
\end{equation}
while 
\begin{equation}\label{case2r}
R_k(I)\le R_{k+1}(I)\quad \text{for}\quad I\ge I_{k+1}.
\end{equation}
Using the notation $\Delta_S$ and $\Delta_I$ introduced in the proof 
of Lemma \ref{l5}, from \eqref{case2} 
we obtain \eqref{Spart} and
$\Delta_S\le 0$,
while \eqref{case2r} leads to
\eqref{kl}, hence the difference \eqref{tildeV} satisfies \eqref{neest'}. Further, relation \eqref{I=I} also holds true,
hence \eqref{lemt'} implies \eqref{neest}.
    Finally, a counterpart of \eqref{cou} reads
        \begin{eqnarray}\nonumber
     &&-\int_{I_*^k}^{I_{k+1}}\frac{f_{k+1}(I_*^{k+1})-f_k(I_*^k)}{f_k(i)}\,di\\
  &&  \nonumber
 \le  \frac{I_*^k-I_{k+1}}{f_k(I_{k+1})}\bigl((I_*^{k+1}-I_*^{k})R_{k+1}(I_*^k) +I_*^k(R_{k+1}(I_*^k)-R_{k}(I_*^k))\bigr)
    \end{eqnarray}
 and due to \eqref{ddt'} and \eqref{lemt'},
 \[
 -\int_{I_*^k}^{I_{k+1}}\frac{f_{k+1}(I_*^{k+1})-f_k(I_*^k)}{f_k(i)}\,di \le 
 \frac{q_0(I_*^k-I_{k+1})^2}{f_k(I_{k+1})}\left( \frac{\rho R_{k+1}(I_*^k)}{R_{k+1}(I_*^{k+1})R_k(I_*^k)}+I_*^k\right).
 \]
   Combining this estimate with \eqref{neest'} and \eqref{neest}, we obtain (cf.\ \eqref{eql5})
   \[
  \tilde \Delta V_k \le q_0(I_*^k-I_{k+1})^2\left(
  \frac{\rho^2 q_0}{R_{k+1}^2(I_*^{k+1})R_k^2(I_*^k) I_*^k}+ 
 \frac{\rho R_{k+1}(I_*^k)}{R_{k+1}(I_*^{k+1})R_k(I_*^k)f_k(I_{k+1})}+\frac{I_*^k}{f_k(I_{k+1})}\right)
   \]
   and adding \eqref{nes} leads to \eqref{VVV}. \hfill $\Box$
   
   \medskip
    
%
    
    
In order to use Lemmas \ref{lmain} and \ref{lmain'} we need to estimate the components $I(t), S(t)$
of the trajectory from below.

  \begin{lemma}\label{lm2}
 For $t\ge t_2$, 
 \begin{equation}\label{estim}
 I(t)\ge i_0:= \rho\,\frac{R_0^{int}-1}{R_0^{int}}\,e^{-\left(\frac{2}{\rho(R_0^{int}-1)}+\frac{1}{R_0^{int}}\right)R_0^{nat}},\qquad
 S(t)\ge s_0:=\frac{e^{-(1+2R_0^{nat})}}{R_0^{nat}}.
 \end{equation}
   \end{lemma}
   
   {\em Proof. } Consider a point $(I_k, S_k)=(I_k,{1}/{R_k(I_k)})$    where  the trajectory satisfies $\dot{I}=0,$ $\dot{S}<0$,
 see Figure \ref{fig2}.
   At this point,
   $$
   V_k(I_k, S_k)=\int_{I_*^k}^{I_k}\left(1-\frac{f_k(I_*^k)}{f_k(i)}\right)di+\int_{S_*^k}^{S_k}\left(1-\frac{S_*^k}{s}\right)ds\le S_k-S_*^k+I_k-I_*^k\le 2.
   $$
 Clearly, the minimal value $s_m^k = \min_{t\in [t_k,t_{k+2}]} S(t)$ of the $S$-component of the trajectory over the time interval 
 $[t_k,t_{k+2}]$ is achieved at a time $t_m^k\in [t_k,t_{k+1}]$, i.e.\ $s_m^k=S(t_m^k)$. Since $V_k$ decreases on the time interval $[t_k,t_{k+1}]$,
from $  V_k(I_k, S_k)\le 2$, it follows that
   $$
   2\ge V_k(I(t_m^k),S(t_m^k))\ge 
   \int_{s_m^k}^{S_*^k}\left(\frac{S_*^k}{s}-1\right)ds=S_*^k\ln{\frac{S_*^k}{s_m^k}}-(S_*^k-s_m^k)\ge S_*^k\ln{\frac{S_*^k}{s_m^k}}-S_*^k.
   $$
   Therefore,
$s_m^k\ge S_*^k e^{-\bigl(1+\frac{2}{S_*^k}\bigr)}$ and due to $S_*^k\ge 1/R_0^{nat}$,
 $$
 s_m^k=\min_{t\in [t_k,t_{k+2}]} S(t) \ge \frac{e^{-(1+2R_0^{nat})}}{R_0^{nat}}.
 $$

  On the other hand, the minimal value $i_m^k = \min_{t\in [t_k,t_{k+2}]} I(t)$ of the $I$-component of the trajectory over the time interval 
 $[t_k,t_{k+2}]$ is achieved at the point $(I_{k+1},S_{k+1})$, i.e.\ $i_m^k=I_{k+1}$. Again,
 since $V_k$ decreases on the time interval $[t_k,t_{k+1}]$,
from $  V_k(I_k, S_k)\le 2$, it follows that
   \begin{eqnarray*}
   2\ge V_k(I_{k+1},S_{k+1}) &\ge& \int_{I_{k+1}}^{I_*^k}\left(\frac{f_k(I_*^k)}{f_k(i)}-1\right)di=f_k(I_*^k)\int_{I_{k+1}}^{I_*^k}\frac{di}{iR_k(i)}-(I_*^k-I_{k+1})
   \\
   &\ge& \frac{f_k(I_*^k)}{R_0^{nat}}\int_{I_{k+1}}^{I_*^k}\frac{di}{i}-I_*^k=\frac{f_k(I_*^k)}{R_0^{nat}}\ln{\frac{I_*^k}{I_{k+1}}}-I_*^k,
   \end{eqnarray*}
hence
  $$ I_{k+1}\ge I_*^ke^{-(2+I_*^k)\frac{R_0^{nat}}{f_k(I_*^k)}}=I_*^ke^{-\left(\frac{2}{I_*^k R_k(I_*^k)}+\frac{1}{R_k(I_*^k)}\right)R_0^{nat}}.$$\\
Therefore, using $I_*^k=\rho(1-S_*^k)$ and $ S_*^k=1/R_k(I_*^k)\le 1/R_0^{int}$, we obtain
\[
      I_{k+1}=\min_{t\in [t_k,t_{k+2}]} I(t)
      \ge 
      \rho(1-S_*^k)e^{-\left(\frac{2}{\rho(1-S_*^k)}+1\right)S_*^k R_0^{nat}}
   \ge\rho\,\frac{R_0^{int}-1}{R_0^{int}}\,e^{-\left(\frac{2}{\rho(R_0^{int}-1)}+\frac{1}{R_0^{int}}\right)R_0^{nat}},
\]
 which completes the proof. \hfill $\Box$
 
    \medskip
    We can now prove the convergence of the trajectories to the equilibrium set.
    
        \begin{lemma}\label{conv1}
  Every trajectory of  system \eqref{r0model}, \eqref{pre} with  initial values $I(0)>0, S(0)>0$
  converges to the set of endemic equilibrium states.
        \end{lemma}
    
    
  {\em Proof.} From Lemmas \ref{lmain}\,--\,\ref{lm2}, it follows that
       the sequence $v_k= V_k(I_k, S_k)$ satisfies
\begin{equation}\label{1}
 v_{k+1}-v_k\le -\kappa(I_*^k-I_{k+1})^2, \qquad k=0,1,2,\ldots
 \end{equation}
%
with
\begin{equation}\label{kappa}
\kappa:=\frac{\rho}{4}\sqrt{\frac{\epsilon_0 s_0}{2R_0^{nat}}}
    -\frac{q_0}{i_0 R_0^{int}}\left(
  \frac{\rho^2 q_0}{(R_0^{int})^3}+ 
 \frac{\rho R_0^{nat}}{(R_0^{int})^2}+1\right)>0,
\end{equation}
 independent of $k$. 
  
Since $f(I)=R(I)I\le R_0^{nat}$ for $0\le I\le 1$, Lemma \ref{l33} implies that
    $$\int_{I_*^k}^{I}\frac{f_k(i)-f_k(I_*^k)}{f_k(i)}di\ge \frac{\epsilon_0}{R_0^{nat}}\int_{I_*^k}^{I}(i-I_*^k)di=\frac{\epsilon_0}{2R_0^{nat}}(I-I_*^k)^2.$$
 On the other hand,
 $$f'_k(I)=(I R_k(I))'=R'_k(I)I+R_k(I)\le R_k(I)\le R^{nat}_0,\quad f_k(I)= I R_k(I)\ge i_0 R_0^{int},$$
where $i_0$ is defined in \eqref{estim}, hence
   $$
   \int_{I_*^k}^{I}\frac{f_k(i)-f_k(I_*^k)}{f_k(i)}di 
    \le 
\frac{R_0^{nat}}{i_0 R_0^{int}}\int_{I_*^k}^{I}(i-I_*^k)di
=\frac{R_0^{nat}}{2i_0 R_0^{int}}(I-I_*^k)^2,
    $$
so we have 
   \begin{equation}\label{a}
 \frac{\epsilon_0}{2R_0^{nat}}(I-I_*^k)^2\le   \int_{I_*^k}^{I}\frac{f_k(i)-f_k(I_*^k)}{f_k(i)}di 
    \le 
\frac{R_0^{nat}}{2i_0 R_0^{int}}(I-I_*^k)^2.
    \end{equation}
    Similarly, due to $0\le S\le 1$,
\begin{equation}\label{b}
 \frac{(S-S_*^k)^2}{2}\le    \int_{S_*^k}^{S}\left(1-\frac{S_*^k}{s}\right)ds\le  \frac{(S-S_*^k)^2}{2s_0},
     \end{equation}
 with $s_0$ defined in \eqref{estim}.    
     Adding \eqref{a} and \eqref{b}, we obtain
          \begin{equation} \label{cc}
     \frac{\epsilon_0(I-I_*^k)^2}{2R_0^{nat}}+\frac{(S-S_*^k)^2}{2}\le  V_k(I,S)\le
     \frac{R_0^{nat}(I-I_*^k)^2}{2i_0 R_0^{int}}+\frac{(S-S_*^k)^2}{2s_0}.
     \end{equation}
    
Relation \eqref{1} implies
\[ 
v_n =v_1 +\sum_{k=1}^{n-1} (v_{k+1}-v_{k})\le
v_1 - \kappa \sum_{k=1}^{n-1} (I_*^k-I_{k+1})^2,
\]
and due to the non-negativity of each Lyapunov function $V_n$, 
    $$\sum_{k=1}^{n-1}(I_*^k-I_{k+1})^2\le \frac{1}{\kappa}v_1$$ for any $n$, hence $$\sum_{k=1}^{\infty}(I_*^k-I_{k+1})^2<\infty.$$\\
    Therefore $I_{k+1}-I_*^k\xrightarrow[]{} 0 $
    (as $ k\xrightarrow[]{} \infty$).
%
    Due to \eqref{l2} and \eqref{lemt'}, this 
    implies 
    $I_k-I_*^k\xrightarrow[]{} 0$, 
and from $S_*^k={1}/{R_k(I_*^k)}$, $S_k={1}/{R_k(I_k)}$ it follows that $S_k-S_*^k\xrightarrow[]{} 0$.
Hence, \eqref{cc} implies
$v_k=V_k(I_k,S_k)\to 0$.
Since $V_k$ is nonnegative and decreases along the trajectory, we conclude that
   \begin{equation}\label{c}
V_k(I_k,S_k)=    \max_{t\in[t_k,  t_{k+1}]} 
   V_k(I(t),S(t))\xrightarrow[]{} 0 \quad \text{as} \quad k\xrightarrow[]{} \infty,
   \end{equation}
%
 %
     hence 
     \eqref{cc} implies that
    the trajectory converges to the equilibrium set.\hfill $\Box$

 
 \medskip
It remains to show that every trajectory converges to an equilibrium.
The following lemma completes the proof of the theorem.

        \begin{lemma}\label{conv2}
  Every trajectory of  system \eqref{r0model}, \eqref{pre} with positive initial values $I(0), S(0)$
  converges to an endemic equilibrium state.
        \end{lemma}

{\em Proof.}
From \eqref{cc}, it follows that
\[
         V_k(I_{k+1}, S_{k+1}) \le \frac{R_0^{nat}}{2i_0 R_0^{int}}(I_{k+1}-I_*^k)^2+\frac{1}{2s_0}\left(\frac{1}{R_k(I_{k+1})}-\frac{1}{R_k(I_*^k)}\right)^2,
 \]        
 where $R_k(I)\ge R_0^{int}$ and $|R_k(I_{k+1})-R_k(I_*^{k})|\le q_0 |I_{k+1}-I_*^k|$ due to \eqref{LipP*}, hence        
 \begin{equation}       \label{VV5}
 V_k(I_{k+1}, S_{k+1})   \le a(I_{k+1}-I_*^k)^2, \qquad a:=\frac{R_0^{nat}}{2i_0 R_0^{int}}+\frac{q_0}{2s_0(R_0^{int})^2}.
\end{equation}
    Lemma 
    \ref{l5} implies that
    $$
    \tilde \Delta V_k \le b  (I_{k+1}-I_*^k)^2, \qquad b:=
\frac{q_0 }{R_0^{int}}  \left(\frac{\rho^2 q_0}{i_0 R_0^{int} (R_0^{int})^2 }+\frac{\rho}{i_0 R_0^{int}}+1\right),
$$
hence,
using \eqref{VV5},
\begin{equation}\label{es1}
V_{k+1}(I_{k+1},S_{k+1})= V_k(I_{k+1}, S_{k+1})  + \tilde \Delta V_k \le (a+b)(I_{k+1}-I_*^k)^2.
\end{equation}

If $V_k(I_k,S_k)\ge (a+b+\kappa)(I_{k+1}-I_*^k)^2$  with $\kappa$ defined by \eqref{kappa}, then \eqref{es1} implies
\begin{equation}\label{es}
v_{k+1}=V_{k+1}(I_{k+1},S_{k+1})\le p\, V_k(I_k,S_k)= p v_k, \qquad p:=\frac{a+b}{a+b+\kappa}<1.
\end{equation}
On the other hand, if $V_k(I_k,S_k)\le (a+b+\kappa)(I_{k+1}-I_*^k)^2,$ then \eqref{1} implies \eqref{es}.
    Therefore, in either case,
    \begin{equation}\label{VVVV}
    v_k \le v_0 p^k.
    \end{equation}
    Since $V_k$  decreases along the trajectory, $V_{k}(I_{k+1},S_{k+1})\le V_k(I_k,S_k)=v_k$.
    But according to \eqref{cc}, 
      \begin{equation}\label{ll}
     V_k(I_{k+1},S_{k+1})\ge \frac{\epsilon_0}{2R_0^{nat}} (I_{k+1}-I_*^k)^2,
     \end{equation}
  hence
     $$
     \frac{\epsilon_0}{2R_0^{nat}}(I_{k+1}-I_{*}^k)^2\le 
     v_k\le v_0 p^k.
     $$
 Further,    
    according to 
    \eqref{lemt}, \eqref{lemt'}, 
   $$|I_*^{k+1}-I_*^k|\le \frac{\rho q_0}{(R_0^{int})^2} |I_{k+1}-I_*^k|,$$ 
  consequently
   $$|I_*^{k+1}-I_*^k|\le \frac{\rho q_0}{(R_0^{int})^2} \sqrt{\frac{2v_0 R_0^{nat}}{\epsilon_0}}\, p^{\frac{k}{2}}.$$
   We see that $I_*^k$ is a Cauchy sequence, hence it converges to a limit value $I_*$.
   Since $I_{k+1}-I_*^k\xrightarrow[]{} 0$ as $k\xrightarrow[]{} \infty$, we conclude that $I(t)\to I_*$ as $t\to\infty$,  which implies that each of the functions $R_0(t)$ and $S(t)$ also has a limit as $t\to \infty$. Denoting the limit of $S(t)$ by $S_*$, the trajectory converges to the equilibrium $(I_*,S_*)$. \hfill $\Box$
   
   

\end{document}